\documentclass[a4paper,11pt]{article}
\pdfoutput=1
\usepackage[utf8]{inputenc}
\usepackage[english]{babel}
\usepackage[T1]{fontenc}
\usepackage[top=3.5cm, bottom=3.5cm, left=3cm, right=3cm]{geometry}
\usepackage{color}    
\usepackage{graphicx} 
\usepackage{listings}
\usepackage{tikz}
\usetikzlibrary{math}
\usepackage{fancyhdr} 
\usepackage{lastpage} 
\usepackage{amsmath,amssymb,mleftright} 
\usepackage{listings} 
\usepackage{parcolumns} 
\usepackage{array,tabu,booktabs} 
\usepackage{multirow} 
\usepackage{pgfplots,wrapfig}
\usepackage{titling}  
\usepackage{lmodern}  
\usepackage{verbatim,verbatimbox} 
\usepackage{fancyvrb} 
\usepackage{enumitem} 
\usepackage{hyperref} 
\usepackage{float}
\usepackage{caption}
\usepackage{subcaption}
\usepackage{animate}
\usepackage[final]{pdfpages}
\usepackage{diagbox}
\usepackage{hhline}
\usepackage{titlesec}
\usepackage{algorithm}
\usepackage{algpseudocode}
\usepackage{bigints}
\usepackage[makeroom]{cancel}
\usepackage{xparse}
\usepackage{eufrak}
\usepackage{setspace}
\usepackage{mathtools}
\titlespacing*{\section}
{0pt}{5.5ex plus 1ex minus .2ex}{4.3ex plus .2ex}
\titlespacing*{\subsection}
{0pt}{5.5ex plus 1ex minus .2ex}{4.3ex plus .2ex}

\usepackage{todonotes} 

\usepackage{titlesec}
\usepackage{rotating}

\titleformat{\paragraph}
{\normalfont\normalsize\bfseries}{\theparagraph}{1em}{}
\titlespacing*{\paragraph}
{0pt}{3.25ex plus 1ex minus .2ex}{1.5ex plus .2ex}

\makeatletter
\newcommand*\bigcdot{\mathpalette\bigcdot@{.5}}
\newcommand*\bigcdot@[2]{\mathbin{\vcenter{\hbox{\scalebox{#2}{$\m@th#1\bullet$}}}}}
\makeatother

\makeatletter
\def\BState{\State\hskip-\ALG@thistlm}
\makeatother
\definecolor{tempcolor}{rgb}{0.7773 0.0820 0.5195}

\newcommand{\lb}{\left[}
\newcommand{\rb}{\right]}

\newcommand{\R}{\mathbb{R}}

\newcommand{\matlab}{\textsc{Matlab}}
\newcommand{\python}{\textsc{Python}}
\newcommand{\fenics}{\textsc{FEniCS}}

\newcommand{\norm}[1]{\ensuremath{\lVert  #1 \rVert}}
\newcommand{\abs}[1]{\ensuremath{\left\vert  #1\right\vert}}

\newcommand{\p}{\partial}

\usepackage{xcolor}
\definecolor{dkgreen}{rgb}{0,0.6,0}
\definecolor{dred}{rgb}{0.545,0,0}
\definecolor{dblue}{rgb}{0,0,0.545}
\definecolor{lgrey}{rgb}{0.9,0.9,0.9}
\definecolor{gray}{rgb}{0.4,0.4,0.4}
\definecolor{darkblue}{rgb}{0.0,0.0,0.6}
\definecolor{turquoise}{rgb}{0.2500,0.8750,0.8125}
\definecolor{indigo}{rgb}{0.2930, 0,0.5078}
\definecolor{mag}{rgb}{1, 0,1}
\definecolor{corn}{rgb}{0.3906,0.5820,0.9258}
\definecolor{mvr}{rgb}{0.7773,0.0820,0.5195}
\definecolor{dod}{rgb}{0.11719,0.5625,1}
\lstdefinelanguage{MatlabCostum}{
      backgroundcolor=\color{white},  
      basicstyle=\footnotesize \ttfamily \color{black} \bfseries,   
      breakatwhitespace=false,       
      breaklines=true,               
      captionpos=b,                   
      commentstyle=\color{dkgreen},   
      emph={repmat, ones},			
      keywordstyle=\color{blue},           
      escapeinside={\%*}{*)},                  
      frame=single,                  
      language=Matlab,                
      identifierstyle=\color{black},
      stringstyle=\color{blue},      
      numbers=left,                 
      numbersep=5pt,                  
      numberstyle=\tiny\color{black}, 
      rulecolor=\color{black},        
      showspaces=false,               
      showstringspaces=false,        
      showtabs=false,                
      stepnumber=1,                   
      tabsize=4,
      title=\lstname
}
\providecommand{\keywords}[1]
{
  \small	
  \textbf{\textit{Keywords---}} #1
}

\newtheorem{theorem}{Theorem}[section]

\newtheorem{definition}[theorem]{Definition}
\newtheorem{lemma}[theorem]{Lemma}
\newtheorem{proposition}[theorem]{Proposition}
\newtheorem{remark}[theorem]{Remark}

\newcommand{\beginproof}{\emph{Proof.\hspace{4mm}}}%
\newcommand*{\QEDA}{\hfill\ensuremath{\blacksquare}}%
\allowdisplaybreaks

\usepackage[backend=biber,
            style=alphabetic,
            abbreviate=false,
            dateabbrev=false,
            alldates=long]{biblatex}

\begin{document}

\title{Jacobian of solutions to the conductivity equation in limited view}

\author{Mikko Salo, Hjørdis Schlüter}
\maketitle


\begin{abstract}
    The aim of hybrid inverse problems such as Acousto-Electric Tomography or Current Density Imaging is the reconstruction of the electrical conductivity in a domain that can only be accessed from its exterior. In the inversion procedure, the solutions to the conductivity equation play a central role. In particular, it is important that the Jacobian of the solutions is non-vanishing. In the present paper we address a two-dimensional limited view setting, where only a part of the boundary of the domain can be controlled by a non-zero Dirichlet condition, while on the remaining boundary there is a zero Dirichlet condition. For this setting, we propose sufficient conditions on the boundary functions so that the Jacobian of the corresponding solutions is non-vanishing. In that regard we allow for discontinuous boundary functions, which requires the use of solutions in weighted Sobolev spaces. We implement the procedure of reconstructing a conductivity from power density data numerically and investigate how this limited view setting affects the Jacobian and the quality of the reconstructions.
\end{abstract}

\keywords{acousto-electric tomography, current density imaging, hybrid inverse problems, coupled physics imaging, non-vanishing Jacobian, conductivity equation}

\section{Introduction}
In certain imaging applications it is important to know whether solutions $u_1$ and $u_2$ to the conductivity equation 
\begin{equation*}
    \begin{cases}
        -\mathrm{div}(\sigma \nabla u_i)=0 & \text{in }\Omega,\\
        u_i=g_i & \text{on }\partial \Omega,
    \end{cases}
\end{equation*}
satisfy the following non-vanishing Jacobian condition:
\begin{equation}
    \mathrm{det}[\nabla u_1(x) \, \nabla u_2(x)]\neq 0, \quad \text{for } x \in \Omega. \label{eq:nonvan}
\end{equation}
Here $\Omega \subset \mathbb{R}^2$ is a bounded Lipschitz domain and $\sigma \in L^{\infty}(\Omega, \mathbb{R}^{2 \times 2})$ is an anisotropic conductivity. This question arises in Acousto-Electric Tomography that aims at reconstructing the unknown interior conductivity $\sigma$ from internal data composed of power density measurements \cite{Zhang:Wang:2004, ammari2008a}. Similar questions appear in other imaging methods including Current Density Imaging \cite{WidlakScherzer, Bal2013_survey, LiSchotlandYangZhong} and Magnetic Resonance Electric Impedance Tomography \cite{Seo_2005,SeoWoo11} that aim at reconstructing the conductivity from current density measurements. These questions are relevant in any dimension $n \geq 2$, but in this article, we will restrict our attention only on the two-dimensional case.

The reconstruction procedure in Acousto-Electric Tomography is characterized by two steps: First reconstructing interior power density data $H_{ij}=\sigma \nabla u_i \cdot \nabla u_j$ from combined information from boundary measurements and perturbations by acoustic waves, and secondly reconstructing $\sigma$ from the power density matrix $\mathbf{H}\in \mathbb{R}^{2 \times 2}$. The non-vanishing Jacobian condition \eqref{eq:nonvan} is essential for the second step in the reconstruction procedure, as it requires inverting the matrix $\mathbf{H}$. \par

The question whether one can find conditions on the boundary functions $g_1$ and $g_2$ so that the non-vanishing Jacobian condition \eqref{eq:nonvan} is satisfied dates back to Radó in the 1920s. For the constant coefficient case $\sigma=\mathbf{I}_2$ an answer to this question is formulated in the Radó-Kneser-Choquet theorem \cite{Rado26,Kneser26,Choquet45}. This result was generalized to non-constant coefficients in \cite{Alessandrini86,Alessandrini87,AlessandriniMagnanini94,AlessandriniNesi01,BaumanMariniNesi01,alessandriniNesi15}. For instance, \cite{BaumanMariniNesi01} require that $g=(g_1,g_2)$ is a $C^1$ diffeomorphism and maps $\partial \Omega$ onto the boundary of a convex domain for the condition \eqref{eq:nonvan} to hold. A discussion of  results of this type is given in \cite{AlbertiCapdeboscq18} (see also \cite{AlbertiArxiv} for recent work on random boundary data). \par

In this paper, we address the same question in a limited view setting that is characterized by a non-empty closed part of the boundary, $\Gamma \subset \partial \Omega$, which we can control, while on the rest of the boundary the potentials $u_1$ and $u_2$ are vanishing:
\begin{equation}\label{bvpui1}
    \begin{cases}
        -\mathrm{div}(\sigma \nabla u_i)=0 & \text{in }\Omega,\\
        u_i=f_i & \text{on }\Gamma ,\\
        u_i = 0 & \text{on }\partial \Omega \backslash \Gamma.
    \end{cases}
\end{equation}
The Radó-Kneser-Choquet type results mentioned above cannot be applied directly in limited view, as $g=(u_1,u_2)\vert_{\partial \Omega}$ is not injective. However, we show that the arguments for proving such results can be adapted to the limited view setting, and we formulate sufficient conditions under which the corresponding Jacobian is non-vanishing. We also allow the boundary functions to be discontinuous (e.g.\ piecewise smooth), which seems natural in this setting and requires the use of weighted Sobolev spaces. For $H^{1/2}$ boundary functions having bounded variation the result is a consequence of \cite[Theorem 2.7]{AlessandriniMagnanini94}, while the result for discontinuous boundary functions is new. We illustrate by numerical simulations how these conditions can be used to reconstruct an isotropic conductivity from power density data. For the numerical simulations an analytic reconstruction approach is used \cite{monard2012a}.

\section{Main results}
We will consider the conductivity equation $-\mathrm{div}(\sigma \nabla u) = 0$ in $\Omega$, where the conductivity matrix $\sigma$ is symmetric and satisfies for some $\lambda,\Lambda > 0$ the ellipticity condition 
\begin{equation} \label{ellip}
\lambda \abs{\xi}^2 \leq \sigma^{jk}(x) \xi_j \xi_k \leq \Lambda \abs{\xi}^2  \text{for a.e.\ $x \in \Omega$ and all $\xi \in \R^n$}.
\end{equation}
The first result states that the presence of an interior critical point for a nonconstant solution $u$ forces oscillations in its boundary data. The result is known for $H^{1/2}$ boundary data, see \cite[Proposition 6.7]{AlbertiCapdeboscq18}, but we give an extension to the case where the boundary data can be slightly worse that $H^{1/2}$ (e.g.\ piecewise smooth). Boundary data in $H^s$ with $s \leq 1/2$ may be discontinuous, and for such functions, we can use a quasicontinuous representative to talk about their pointwise values \cite[Chapter 6]{AdamsHedberg}. The notation $H^1(\Omega,d^{1-2s})$ for weighted spaces is explained in Section \ref{sec_weighted}.

\begin{proposition}
\label{PropNoCriticalPoints}
Let $\Omega \subset \R^2$ be a simply connected bounded Lipschitz domain and let $\sigma \in C^{0,\alpha}(\overline{\Omega})$ satisfy \eqref{ellip}. There is $\varepsilon > 0$ with the following property: if $f \in H^s(\partial \Omega)$ where $\abs{s-1/2} < \varepsilon$ and if $u \in H^1(\Omega, d^{1-2s})$ is a nonconstant solution of 
\begin{equation*}
\begin{cases}
    -\mathrm{div}(\sigma \nabla u)=0 & \text{in }\Omega,\\
    u = f & \text{on }\partial \Omega,
\end{cases}
\end{equation*}
and if $\nabla u(x_0) = 0$ for some $x_0 \in \Omega$, then there are $x_1, x_2, x_3, x_4 \in \partial \Omega$ that are in counterclockwise order along $\partial \Omega$ such that 
\[
u(x_1) > u(x_0), \quad u(x_2) < u(x_0), \quad u(x_3) > u(x_0), \quad u(x_4) < u(x_0).\]
\end{proposition}
\beginproof We follow the argument in \cite[Proposition 6.7]{AlbertiCapdeboscq18}. By \cite[Theorem 2.3.3]{schulz1990} the interior regularity of $u$ is $C^{1,\alpha}_{loc}(\Omega)$. As $x_0$ is a critical point of $u$ it follows from \cite[Proposition 6.6]{AlbertiCapdeboscq18} that in a neighborhood $U$ of $x_{0}$ the level set $\{x \in U \,:\, u(x)=u(x_{0})\}$ is made of $m+1$ arcs intersecting with equal angles at $x_{0}$ for some $m\geq 1$. We note that by \cite[Proposition 6.5 (i)]{AlbertiCapdeboscq18} the set $\{x \in U \,:\, u(x)>u(x_{0})\}$ is made of $m+1$ connected components that we denote by $U_l^+$:
\begin{equation*}
    \{x \in U \,:\, u(x)>u(x_{0})\}=\bigcup_{l=1}^{m+1} U_l^+.
\end{equation*}
Furthermore, by the same proposition it follows that these connected components alternate with the corresponding connected components $U_l^-$ of the set $\{x \in U \,:\, u(x)<u(x_{0})\}$. We now consider the corresponding sets over the whole domain $\Omega$. The components of $\{ u(x)>u(x_{0}) \}$ are denoted by $\Omega_j^+$ and 
\begin{equation*}
    \{x \in \Omega \,:\, u(x)>u(x_{0})\}=\bigcup_{j\in J} \Omega_j^+.
\end{equation*}
Similarly, the components of $\{ u(x)<u(x_{0}) \}$ are denoted by $\Omega_j^-$ and 
\begin{equation*}
    \{x \in \Omega \,:\, u(x)<u(x_{0})\}=\bigcup_{k\in K} \Omega_k^-.
\end{equation*}
Now pick indices $j_1,j_2 \in J$ such that $U_1^+ \subset \Omega_{j_1}^+$ and $U_2^+ \subset \Omega_{j_2}^+$. It follows from theorem \ref{thm_mp} that the weak maximum principle holds for $H^1(\Omega,d^{1-2s})$ solutions (in order to apply theorem \ref{thm_mp} it is required that $s$ satisfies $\abs{s-1/2}<\varepsilon$ for a suitable $\varepsilon>0$). Therefore if $u(x)>u(x_0)$ holds in the interior of the domains $\Omega_{j_1}^+$ and $\Omega_{j_2}^+$, then one must have $u(x_1) > u(x_0)$ and $u(x_3) > u(x_0)$ for some $x_1 \in \partial \Omega_{j_1}^+$ and $x_3 \in \partial \Omega_{j_2}^+$. Since $u(x) = u(x_0)$ for $x \in \partial \Omega_{j_i}^+ \cap \Omega$, we must have $x_{1}\in \partial \Omega \cap \overline{\Omega}_{j_1}^+$ and $x_{3}\in \partial \Omega \cap \overline{\Omega}_{j_2}^+$. By construction there exist indices $l_1,l_2 \in [1,..,m+1]$ so that $U_{l_1}^-$ is located between $U_1^+$ and $U_2^+$ and $U_{l_2}^-$ is located to the other side of $U_2^+$. We now pick indices $k_1,k_2 \in K$ such that $U_{l_1}^-\subset \Omega_{k_1}^-$ and $U_{l_2}^-\subset \Omega_{k_2}^-$. By the weak maximum principle (theorem \ref{thm_mp}) there then exist points $x_2\in \partial \Omega \cap \overline{\Omega}_{k_1}^-$ and $x_4\in \partial \Omega \cap \overline{\Omega}_{k_2}^-$ such that $u(x_2)<u(x_0)$ and $u(x_4)<u(x_0)$ yielding the desired statement.
\QEDA\\

Let $\gamma: [0,\ell] \to \R^2$ be a $C^1$ curve (we do not require that $\gamma(0) = \gamma(\ell)$). We say that $\gamma$ is regular if $\dot{\gamma}(t) \neq 0$ for all $t \in [0,\ell]$. For a regular curve, we may write a polar coordinate representation for the tangent vector $\dot{\gamma}(t)$ as 
\[
\dot{\gamma}(t) = r(t) e^{i\phi(t)}
\]
where $r(t) = \abs{\dot{\gamma}(t)}$ and $\phi(t)$ are continuous functions in $[0,\ell]$. The function 
\[
\mathrm{arg}(\dot{\gamma}(t)) := \phi(t)
\]
is well defined modulo a constant in $2\pi \mathbb{Z}$. We define 
\[
\mathrm{Ind}(\dot{\gamma}) := \frac{\mathrm{arg}(\dot{\gamma}(\ell)) - \mathrm{arg}(\dot{\gamma}(0))}{2\pi}.
\]
If $\gamma$ is closed, i.e.\ $\gamma(0) = \gamma(\ell)$, then $\mathrm{Ind}(\dot{\gamma})$ is the winding number of the curve $\dot{\gamma}(t)$ (also called the rotation index of $\gamma(t)$), which is an integer. If $\gamma$ is not closed but $\mathrm{arg}(\dot{\gamma}(t))$ is monotone (i.e.\ nondecreasing or nonincreasing), we may still interpret $\mathrm{Ind}(\dot{\gamma})$ as the winding number of $\dot{\gamma}(t)$, and this is then a real number.

We now give sufficient conditions on a pair of boundary data vanishing outside an arc $\Gamma$ such that the corresponding solutions $u_1, u_2$ satisfy $\mathrm{det} [\nabla u_1(x) \, \nabla u_2(x)]\neq 0$ everywhere in $\Omega$. Condition (a) below is related to the case where $u_i|_{\partial \Omega}$ are continuous and condition (b) allows discontinuous boundary data. Note that part (a) can be seen as a consequence of \cite[Theorem 2.7]{AlessandriniMagnanini94}, while for part (b) Proposition \ref{PropNoCriticalPoints} is needed.

\begin{theorem}\label{thm:Main}
Let $\Omega \subset \R^2$ be a bounded simply connected domain with $C^1$ boundary curve $\eta: [0,2\pi] \to \partial \Omega$, and let $\sigma \in C^{0,\alpha}(\overline{\Omega};\R^{2\times 2})$ satisfy \eqref{ellip}. Let $\Gamma = \eta([0,\ell])$ be a closed arc in $\partial \Omega$. Let $f_1, f_2 \in C^1(\Gamma)$ be linearly independent‚ and assume that $u_i$ is the unique solution of 
\begin{equation}\label{bvpui}
    \begin{cases}
        -\mathrm{div}(\sigma \nabla u_i)=0 & \text{in }\Omega,\\
        u_i=f_i & \text{on }\Gamma ,\\
        u_i = 0 & \text{on }\partial \Omega \backslash \Gamma.
    \end{cases}
\end{equation}
Assume that the curve $\gamma: [0,\ell] \to \R^2$, $\gamma(t) = (f_1(\eta(t)), f_2(\eta(t)))$ is regular, $\mathrm{arg}(\dot{\gamma}(t))$ is monotone, and that one of the following holds:
\begin{enumerate}
    \item[(a)] $u_i|_{\partial \Omega}$ are continuous, and $\abs{\mathrm{Ind}(\dot{\gamma})} \leq 1$; or 
    \item[(b)] $u_i|_{\partial \Omega}$ are continuous at $\eta(0)$, and $\abs{\mathrm{Ind}(\dot{\gamma})} \leq 1/2$.
\end{enumerate}
Then $\mathrm{det} [\nabla u_1(x) \, \nabla u_2(x)]\neq 0$ for all $x \in \Omega$.

\end{theorem}
\beginproof
Assume that (a) or (b) holds, but one has $\mathrm{det} [\nabla u_1(x_0) \, \nabla u_2(x_0)] =  0$ for some $x_0 \in \Omega$. Then there is a vector $\vec{\alpha} = (\alpha_1, \alpha_2) \in \R^2 \setminus \{0\}$ such that the function 
\[
u = \alpha_1 u_1 + \alpha_2 u_2
\]
satisfies $\nabla u(x_0) = 0$. Note that if $u$ is a constant, then by the boundary condition one has $u \equiv 0$ and hence $f_1$ and $f_2$ would be linearly dependent. Thus, we may assume that $u$ is nonconstant. Note that $u_i|_{\partial \Omega}$ are piecewise $C^1$ and hence they are also in $H^s(\partial \Omega)$ for any $s < 1/2$. This implies that $u_i \in H^1(\Omega, d^{1-2s})$ by Theorem \ref{thm_weighted_trace2}. Then by Proposition \ref{PropNoCriticalPoints} there exist distinct points $x_1, x_2, x_3, x_4 \in \partial \Omega$ that are in counterclockwise order along $\partial \Omega$ such that 
\[
u(x_1) > u(x_0), \quad u(x_2) < u(x_0), \quad u(x_3) > u(x_0), \quad u(x_4) < u(x_0).
\]
Consider the function $g: [0,\ell] \to \R$ given by 
\begin{equation}\label{g_fun}
g(t) := u(\eta(t)) = \vec{\alpha} \cdot \gamma(t).
\end{equation}
Extend $g$ by zero to $[0,2\pi)$. Writing $x_j = \eta(t_j)$ where $t_j \in [0,2\pi)$, we have 
\begin{equation} \label{g_contra}
g(t_1) > u(x_0), \quad g(t_2) < u(x_0), \quad g(t_3) > u(x_0), \quad g(t_4) < u(x_0).
\end{equation}
We may assume that $t_1 < t_2 < t_3 < t_4$ (possibly after a cyclic permutation of the indices and after changing $g$ to $-g$).

We now assume that (a) holds, and want to derive a contradiction with \eqref{g_contra}. The function $g$ is $C^1$ on $[0,\ell]$ and its derivative satisfies 
\begin{equation} \label{gprime_formula}
g'(t) = \vec{\alpha} \cdot \dot{\gamma}(t).
\end{equation}
Since $\mathrm{arg}(\dot{\gamma}(t))$ is monotone and $\abs{\mathrm{Ind}(\dot{\gamma}(t))} \leq 1$, it follows that $g'(t)$ either has at most two zeros in $[0,\ell]$, or has three zeros two of which are at $t=0$ and $t=\ell$. Note that if the argument is not strictly monotone, we make the interpretation that some of these zeros of $g'$ could be intervals. We also note that by \eqref{gprime_formula} and monotonicity of $\mathrm{arg}(\dot{\gamma}(t))$, $g'$ changes sign after each of these zeros. Now suppose that $u(x_0) \geq 0$. Using the assumption that $u|_{\partial \Omega}$ is continuous, we have $g(0) = 0$, and then  \eqref{g_contra} implies that $g'$ is positive somewhere in $(0,t_1)$, negative somewhere in $(t_1,t_2)$, positive somewhere in $(t_2,t_3)$, and negative somewhere in $(t_3,\ell)$. On the other hand, if $u(x_0) < 0$, we use the fact that $g(\ell) = 0$ to obtain similarly that $g'$ is negative somewhere in $(t_1,t_2)$, positive somewhere in $(t_2,t_3)$, negative somewhere in $(t_3,t_4)$, and positive somewhere in $(t_4,\ell)$. In both cases $g'$ has at least three zeros in $(0,\ell)$. Moreover, before the first such zero, after the last zero, and between each subsequent pair of these zeros there are points where $g'$ is nonzero. This is a contradiction.

Assume that (b) holds. One has the formula \eqref{gprime_formula} for $g'(t)$ on $[0,\ell]$. Since $\mathrm{arg}(\dot{\gamma}(t))$ is monotone and $\abs{\mathrm{Ind}(\dot{\gamma}(t))} \leq 1/2$, $g'(t)$ either has at most one zero (that could be an interval) in $[0,\ell]$, or has two zeros (that could be intervals) which are at $t=0$ and $t=\ell$. By the assumption that $u_i\vert_{\partial \Omega}$ is continuous at $t=0$ it follows that $g(0)=0$, while there may be a discontinuity at $t=\ell$. If one has $u(x_0) \geq 0$, it follows from \eqref{g_contra} that $g'$ is positive somewhere in $(0,t_1)$, negative somewhere in $(t_1,t_2)$, and positive somewhere in $(t_2,t_3)$. On the other hand if $u(x_0) < 0$, from \eqref{g_contra} one sees that $t_4 \in (0,\ell]$ and hence $g'$ is negative somewhere in $(t_1, t_2)$, positive somewhere in $(t_2,t_3)$, and negative somewhere in  $(t_3,t_4)$. In both cases $g'$ has at least two zeros in $(0,\ell)$ and before, between, and after these zeros there are points where $g'$ is nonzero. This is a contradiction.
\QEDA

\begin{remark}\label{rem:Discont}
In the setting of theorem \ref{thm:Main}, let $\Gamma=\eta([0,\ell])$ and $\Gamma^d=\eta([0,\ell^d])$ with $\ell^d < \ell$. Boundary functions $f_1,f_2\in C^1(\Gamma)$ that satisfy the assumptions and condition (a) in theorem \ref{thm:Main} can also be used to generate boundary functions $f_1^d, f_2^d \in C^1(\Gamma^d)$ whose zero extensions are discontinuous.  Define $f_i^d$ as
\[
    f_i^d(\eta(t))= f_i(\eta(t)) \chi_{[0,\ell^d]}(t).
\]
This yields solutions $u_i^d$ that satisfy $\mathrm{det} [\nabla u_1^d(x) \, \nabla u_2^d(x)]\neq 0$. This allows for boundary functions that are not captured in condition (b) in theorem \ref{thm:Main}, as in this case it is possible that $1/2<\lvert \mathrm{Ind}(\dot{\gamma^d}) \rvert$, where $\gamma^d(t) = (f_1^d(\eta(t)), f_2^d(\eta(t)))$.
\end{remark}

\beginproof
Assume that one has $\det[\nabla u_1^d(x_0) \, \nabla u_2^d(x_0)]=0$ for some $x_0 \in \Omega$. Then there is a vector $\vec{\alpha} = (\alpha_1, \alpha_2) \in \R^2 \setminus \{0\}$ such that the function 
\[
u^d = \alpha_1 u_1^d + \alpha_2 u_2^d
\]
satisfies $\nabla u^d(x_0) = 0$. 
As $u_i^d \vert_{\partial \Omega}$ is piecewise $C^1$ it follows by the analysis in the proof of theorem \ref{thm:Main} that there exist distinct points $x_1, x_2, x_3, x_4 \in \partial \Omega$ that are in counterclockwise order along $\partial \Omega$ such that
\[
u^d(x_1) > u^d(x_0), \quad u^d(x_2) < u^d(x_0), \quad u^d(x_3) > u^d(x_0), \quad u^d(x_4) < u^d(x_0).
\]
We then consider the function $g^d: [0,\ell] \rightarrow \mathbb{R}$,
\begin{equation*}
    g^d(t) := u^d(\eta(t)) = (\vec{\alpha} \cdot \gamma(t)) \chi_{[0,\ell^d]}(t)
\end{equation*}
where $\gamma(t) = (f_1(\eta(t)), f_2(\eta(t))$ for $t \in [0,\ell]$ as before. We extend $g^d$ by zero to $[0,2\pi)$. Writing $x_j = \eta(t_j)$ where $t_j \in [0,2\pi)$, we have 
\begin{equation} \label{g_contra2}
g^d(t_1) > u^d(x_0), \quad g^d(t_2) < u^d(x_0), \quad g^d(t_3) > u^d(x_0), \quad g^d(t_4) < u^d(x_0).
\end{equation}
Furthermore, we consider the function $g: [0,\ell] \rightarrow \mathbb{R}$ for the same vector $\vec{\alpha}$:
\begin{equation*}
    g(t) := \vec{\alpha} \cdot \gamma(t).
\end{equation*}

Since $\mathrm{\arg}(\dot{\gamma}(t))$ is monotone and $\abs{\mathrm{Ind}(\dot{\gamma})} \leq 1$, it follows that $g'(t)$ has at most two zeros in $(0,\ell)$, or three zeros two of which are at $t=0$ and $t=\ell$. (Again these zeros could be intervals.) Since $u_j|_{\partial \Omega}$ are continuous, we have $f_j(\eta(0)) = f_j(\eta(\ell)) = 0$ and thus  $g(0)=g(\ell)=0$. Suppose that $g'(t)$ has exactly two zeros in $(0,\ell)$. Then since $g'$ must change sign after each zero, there exist two intervals $(0,t_{i})$ and $(t_{i},\ell)$ so that either
\[
g(t) \geq 0 \quad \text{for }t \in (0,t_{i}) \quad \text{and} \quad g(t) \leq 0 \quad \text{for }t \in (t_{i},\ell),
\]
or 
\[
g(t) \leq 0 \quad \text{for }t \in (0,t_{i}) \quad \text{and} \quad g(t) \geq 0 \quad \text{for }t \in (t_{i},\ell).
\]
On the other hand if $g'(t)$ has at most one zero, or three zeros two of which are at $t=0$ and $t=\ell$, then either 
\[
g(t) \geq 0 \quad \text{for } t \in (0,\ell) \quad \text{or} \quad g(t) \leq 0 \quad \text{for } t \in (0,\ell).
\]
The behavior of $g$ translates to $g^d$, as $g^d$ is the restriction of $g$ to the interval $[0,\ell^d]$ with $\ell^d < \ell$. It follows that there are at most two intervals for which $g^d$ is nonnegative and nonpositive respectively, and additionally $(g^d)'$ must change sign after each of its zeros. This is in contradiction with \eqref{g_contra2} as no matter if $u^d(x_0)\geq 0$ or $u^d(x_0) \leq 0$ it is not possible for $g^d$ to have two points for which $g^d(t) \geq u(x_0)$ and two points for which $g^d(t) \leq u(x_0)$ in alternating order.  
\QEDA

\begin{remark}
For boundary functions $f_1$ and $f_2$ that satisfy one of the conditions in theorem \ref{thm:Main} it is determined by the order of the functions whether $\mathrm{det}[\nabla u_1(x) \, \nabla u_2(x)]$ is positive or negative for all $x \in \Omega$.
\end{remark}

\begin{remark} \label{rem:vio}
Let $\Omega$ be a $C^{1,\alpha}$ domain and $f_1,f_2 \in C^{1,\alpha}(\overline{\Omega})$, then $u_i \in C^{1,\alpha}$ on $\overline{\Omega}$ away from the endpoints of $\Gamma$ \cite[Corollary 8.36]{Gilbarg:Trudinger:2001}. Due to the limited view setting $\mathrm{det}[\nabla u_1(x) \, \nabla u_2(x)]$ is zero for $x\in \partial \Omega \setminus \Gamma$. Therefore
\begin{equation*}
    \inf_{x \in \Omega} \mathrm{det}[\nabla u_1(x) \, \nabla u_2(x)] = 0.
\end{equation*}
\end{remark}

\beginproof
We decompose $\nabla u_i$ into two parts with contribution from the unit outward normal $\nu$ and the tangent vector $\omega= \mathcal{J} \nu$ with $\mathcal{J}=\begin{bmatrix} 0 & -1\\ 1 & 0\end{bmatrix}$:
\begin{equation*}
    \nabla u_i = (\nabla u_i \cdot \nu)\nu + (\nabla u_i \cdot \omega)\omega.
\end{equation*}
As $u_i\vert_{\partial \Omega \setminus \Gamma}=0$ by the boundary value problem \eqref{bvpui} it follows that there is no contribution in $\omega$-direction. Therefore both $\nabla u_1\vert_{\partial \Omega \setminus \Gamma}$ and $\nabla u_2\vert_{\partial \Omega \setminus \Gamma}$ are parallel to the unit normal $\nu$ so that $\mathrm{det}[\nabla u_1 \vert_{\partial \Omega \setminus \Gamma} \, \nabla u_2\vert_{\partial \Omega \setminus \Gamma}]=0$. 
\QEDA

\section{Dirichlet problem in weighted spaces} \label{sec_weighted}

Let $\Omega \subset \R^n$ be a bounded open set with Lipschitz boundary, and consider the operator 
\[
Lu = -\partial_j(\sigma^{jk} \p_k u) + cu
\]
where $\sigma^{jk}, c \in L^{\infty}(\Omega)$, $\sigma^{jk} = \sigma^{kj}$, and $(\sigma^{jk})$ is uniformly elliptic in the sense that for some $\lambda, \Lambda > 0$, 
\begin{equation} \label{ajk_elliptic}
\lambda \abs{\xi}^2 \leq \sigma^{jk}(x) \xi_j \xi_k \leq \Lambda \abs{\xi}^2  \text{for a.e.\ $x \in \Omega$ and all $\xi \in \R^n$}.
\end{equation}
We wish to consider the Dirichlet problem 
\[
Lu = 0 \text{ in $\Omega$}, \quad u = f \text{ on $\partial \Omega$}
\]
in suitable weighted Sobolev spaces. For general references on weighted Sobolev spaces see \cite[Chapter 3]{Triebel1978} and \cite{Kufner1980}. The following theory in the $L^2$ setting is mostly in \cite{Kufner1980}, but for completeness and possible later use for discontinuous boundary functions we also discuss the $L^p$ theory following \cite{Doyoon08} but with slightly different notation. The results for $p \neq 2$ are not used in the other sections of this article.

\begin{definition}
Let $\Omega \subset \R^n$ be a bounded open set with Lipschitz boundary. Let $1 < p < \infty$ and $\alpha \in \R$, and let $d(x) = \mathrm{dist}(x, \p \Omega)$. Consider the norms 
\begin{align*}
\norm{u}_{L^p(\Omega, d^{\alpha})} &= \lVert u d^{\alpha/p} \rVert_{L^p(\Omega)}, \\
\norm{u}_{W^{1,p}(\Omega,d^{\alpha})} &= \norm{u}_{L^p(\Omega, d^{\alpha})} + \norm{\nabla u}_{L^p(\Omega, d^{\alpha})}.
\end{align*}
Let $W^{1,p}(\Omega, d^{\alpha})$ be the space of all $u \in L^p_{\mathrm{loc}}(\Omega)$ with $\norm{u}_{W^{1,p}(\Omega,d^{\alpha})} < \infty$. We also define $W^{1,p}_0(\Omega, d^{\alpha})$ as the closure of $C^{\infty}_c(\Omega)$ in $W^{1,p}(\Omega, d^{\alpha})$.
\end{definition}

The spaces $W^{1,p}(\Omega, d^{\alpha})$ and $W^{1,p}_0(\Omega, d^{\alpha})$ are Banach spaces, and they are equal when $\alpha \leq -1$ or $\alpha > p-1$ (see \cite[Proposition  9.10]{Kufner1980}). For $\alpha > -1$ the set $C^{\infty}(\overline{\Omega})$ is dense in $W^{1,p}(\Omega, d^{\alpha})$ (see \cite[Remark 7.2]{Kufner1980}). The trace space of $W^{1,p}(\Omega, d^{\alpha})$ can then be identified with a Sobolev space on $\p \Omega$ as follows.

For $1 < p < \infty$ and $0 < s < 1$, let $W^{s,p}(\partial \Omega)$ be the standard Sobolev space on $\p \Omega$ defined via a partition of unity, $C^1$ boundary flattening transformations, and corresponding spaces on $\R^{n-1}$. Part (a) of the following trace theorem is given in \cite[Theorem 2.13]{Doyoon08} (see \cite[Section 3.6.1]{Triebel1978} for the case of $C^{\infty}$ domains), and part (b) follows from \cite[Lemma 2.14, Remark 2.15 and Proposition 2.3]{Doyoon08} and Lemma \ref{lemma_hardy} below.

\begin{theorem} \label{thm_trace}
Let $\Omega \subset \R^n$ be a bounded Lipschitz domain, let $1 < p < \infty$, and let $-1 < \alpha < p-1$.
\begin{enumerate}
\item[(a)]
The trace operator $T: C^{\infty}(\overline{\Omega}) \to C(\p \Omega), \ Tu = u|_{\p \Omega}$ extends as a bounded surjective operator 
\[
T: W^{1,p}(\Omega, d^{\alpha}) \to W^{1 - \frac{1+\alpha}{p},p}(\p \Omega).
\]
Moreover, $T$ has a bounded right inverse $E: W^{1 - \frac{1+\alpha}{p},p}(\p \Omega) \to W^{1,p}(\Omega, d^{\alpha})$.

\item[(b)]
The space $W^{1,p}_0(\Omega, d^{\alpha})$ satisfies 
\begin{align*}
W^{1,p}_0(\Omega, d^{\alpha}) &= \{ u \in W^{1,p}(\Omega, d^{\alpha}) \,:\, Tu = 0\} \\
 &= \{ u \in L^p(\Omega, d^{\alpha-p}) \,:\, \nabla u \in L^p(\Omega, d^{\alpha}) \}.
\end{align*}
The three norms $\norm{\,\cdot\,}_{W^{1,p}(\Omega, d^{\alpha})}$, $\norm{\,\cdot\,}_{L^p(\Omega, d^{\alpha-p})} + \norm{\nabla \,\cdot\,}_{L^p(\Omega, d^{\alpha})}$, and $\norm{\nabla \,\cdot\,}_{L^p(\Omega, d^{\alpha})}$ are equivalent norms on $W^{1,p}_0(\Omega, d^{\alpha})$.
\end{enumerate}
\end{theorem}

The following Hardy inequality is given in \cite[Section 9.1]{Kufner1980}. However, for later purposes we need to make sure that the constant has a controlled dependence on $\alpha$ and hence we repeat the proof.

\begin{lemma} \label{lemma_hardy}
Let $\Omega \subset \R^n$ be a bounded Lipschitz domain, $1 < p < \infty$, and $\alpha \in \R$, $\alpha \neq p-1$. There are $C, C_1 > 0$ only depending on $\Omega, n, p$ such that  
\[
\norm{d^{(\alpha/p)-1} u}_{L^p(\Omega)} \leq C C_1^{\alpha/p} \left( 1 + \frac{1}{\abs{\alpha-p+1}} \right) \norm{d^{\alpha/p} \nabla u}_{L^p(\Omega)}
\]
for any $u \in W^{1,p}_0(\Omega, d^{\alpha})$.
\end{lemma}
\beginproof
We begin with the case of $\R^n_+ = \{ x_n > 0 \}$. Let $u \in C^{\infty}_c(\R^n_+)$. We integrate by parts over $\R^n_+$ and use the H\"older inequality to obtain 
\begin{align*}
    \int x_n^{\alpha-p} u^p \,dx &= \int \p_n\left( \frac{x_n^{\alpha-p+1}}{\alpha-p+1}\right) u^p \,dx = -\frac{p}{\alpha-p+1} \int x_n^{\alpha-p+1-\alpha/p} u^{p-1} x_n^{\alpha/p} \p_n u \,dx \\
    &\leq \frac{p}{\abs{\alpha-p+1}} \norm{x_n^{(\alpha/p)-1} u}_{L^p}^{p-1} \norm{x_n^{\alpha/p} \p_n u}_{L^p}.
\end{align*}
This implies that for any $u \in C^{\infty}_c(\R^n_+)$, one has 
\begin{equation*}
\norm{x_n^{(\alpha/p)-1} u}_{L^p} \leq \frac{p}{\abs{\alpha-p+1}} \norm{x_n^{\alpha/p} \p_n u}_{L^p}.
\end{equation*}
Similarly, if $U = \{ (x',x_n) \in \R^n \,:\, \abs{x'} < r,\ x_n > h(x') \}$ where $r > 0$ and $h: \{ \abs{x'} \leq r\} \to \R$ is a Lipschitz function, and if $u \in C^{\infty}(\overline{U})$ vanishes near $\{x_n = h(x')\}$ and $\{x_n=\infty\}$, the same argument gives that 
\begin{equation} \label{hardy_lipschitzgraph}
\norm{(x_n-h(x'))^{(\alpha/p)-1} u}_{L^p(U)} \leq \frac{p}{\abs{\alpha-p+1}} \norm{(x_n-h(x'))^{\alpha/p} \p_n u}_{L^p(U)}.
\end{equation}

Now if $\Omega$ is a bounded Lipschitz domain, then $\partial \Omega$ can be covered by finitely many balls $B_1, \ldots, B_N$ such that for each $j$, after a rigid motion one has $B_j \cap \Omega = \{ x_n > h_j(x') \} \cap \Omega$ where $h_j$ is a Lipschitz function, and $d(x)$ in $B_j \cap \Omega$ is comparable to $x_n - h_j(x')$ (see \cite[Corollary 4.8]{Kufner1980}). There is also an open set $B_0$ with $\overline{B}_0 \subset \Omega$ so that $\overline{\Omega}$ is covered by $B_0, \ldots, B_N$. Moreover, $B_j$ and $h_j$ only depend on $\Omega$ and not on $p$ and $\alpha$.

Let $u \in C^{\infty}_c(\Omega)$. We can now apply \eqref{hardy_lipschitzgraph} in $B_j \cap \Omega$ for $j=1,\ldots,N$ to obtain that 
\[
\norm{d^{(\alpha/p)-1} u}_{L^p(B_j \cap \Omega)} \leq C \frac{C_1^{\alpha/p}}{\abs{\alpha-p+1}} \norm{d^{\alpha/p} \nabla u}_{L^p(\Omega)}.
\]
In $B_0$, where $d(x)$ is comparable to $1$, we can apply a Poincar\'e inequality as in \cite[Section 9.1]{Kufner1980} and use the above estimates on $B_j \cap \Omega$ to get 
\begin{align*}
\norm{d^{(\alpha/p)-1} u}_{L^p(B_0)} &\leq C C_1^{\frac{\alpha}{p}} \norm{u}_{L^p(B_0)} \leq C C_1^{\alpha/p} (\norm{\nabla u}_{L^p(B_0)} + \norm{u}_{L^p(B_0 \cap (B_1 \cup \ldots \cup B_N)}) \\
 &\leq C C_1^{\alpha/p} \left( 1 + \frac{1}{\abs{\alpha-p+1}} \right) \norm{d^{\alpha/p}\nabla u}_{L^p(\Omega)}.
\end{align*}
The result follows by adding these inequalities and using that $C^{\infty}_c(\Omega)$ is dense in $W^{1,p}_0(\Omega, d^{\alpha})$.
\QEDA \\

The next result, which follows from \cite[Theorem 3.7]{Doyoon08}, states the solvability of the Dirichlet problem in weighted Sobolev spaces when the Dirichlet data is in $W^{s,p}(\p \Omega)$. 

\begin{theorem} \label{thm_weighted_trace}
Let $\Omega \subset \R^n$ be a bounded $C^1$ domain, let $1 < p < \infty$, and let $0 < s < 1$. Assume that $\sigma^{jk}$ and $c$ are Lipschitz continuous in $\overline{\Omega}$ with $(\sigma^{jk})$ satisfying \eqref{ajk_elliptic}, and assume that $c \geq 0$. Given any $f \in W^{s,p}(\p \Omega)$, there is a unique solution $u \in W^{1,p}(\Omega, d^{p(1-s)-1})$ of the problem 
\[
Lu = 0 \text{ in $\Omega$}, \quad u|_{\p \Omega} = f.
\]
One has the estimate 
\[
\norm{u}_{W^{1,p}(\Omega, d^{p(1-s)-1})} \leq C \norm{f}_{W^{s,p}(\p \Omega)}
\]
with $C$ independent of $f$.
\end{theorem}

For $p=2$ we obtain a similar result in weighted $L^2$ spaces $H^1(\Omega, d^{\alpha}) := W^{1,2}(\Omega, d^{\alpha})$ under weaker conditions, but assuming that $s$ is close to $1/2$.

\begin{theorem} \label{thm_weighted_trace2}
Let $\Omega \subset \R^n$ be a bounded Lipschitz domain. Assume that $\sigma^{jk}, c \in L^{\infty}(\Omega)$ with $(\sigma^{jk})$ satisfying \eqref{ajk_elliptic} and $c \geq 0$ a.e.\ in $\Omega$. There is $\varepsilon > 0$ such that whenever $\abs{s-1/2} < \varepsilon$, then for any $f \in H^s(\p \Omega)$ there is a unique solution $u \in H^1(\Omega, d^{1-2s})$ of the problem 
\[
Lu = 0 \text{ in $\Omega$}, \quad u|_{\p \Omega} = f.
\]
One has the estimate 
\[
\norm{u}_{H^1(\Omega, d^{1-2s})} \leq C \norm{f}_{H^s(\p \Omega)}
\]
with $C$ independent of $f$.
\end{theorem}
\beginproof
Note that $\abs{s-1/2} < \varepsilon$ implies $\abs{1-2s} < 2 \varepsilon$. If $\varepsilon$ is chosen small enough, the result follows by combining \cite[Theorem 14.4]{Kufner1980} and the trace theorem (Theorem \ref{thm_trace}) above.
\QEDA \\

The next result gives a weak maximum principle for solutions in $H^1(\Omega, d^{\alpha})$ when $\abs{\alpha}$ is sufficiently small. This smallness condition is analogous to the condition that $s$ is close to $1/2$ in Theorem \ref{thm_weighted_trace2}.

\begin{theorem} \label{thm_mp}
Let $\Omega \subset \R^n$ be a bounded Lipschitz domain. Let $\sigma^{jk}, c \in L^{\infty}(\Omega)$ be such that \eqref{ajk_elliptic} holds and $c \geq 0$ a.e.\ in $\Omega$. There is $\varepsilon > 0$ such that if $\abs{\alpha} \leq \varepsilon$ and $u \in H^1(\Omega, d^{\alpha})$ solves 
\[
-\p_k(\sigma^{jk} \p_j u) + cu = 0 \text{ in $\Omega$}
\]
in the sense of distributions, and if $Tu \leq C$ a.e.\ on $\p \Omega$, then $u \leq C$ in $\Omega$. Similarly, if $Tu \geq C$ a.e.\ on $\p \Omega$, then $u \geq C$ in $\Omega$.
\end{theorem}

The proof uses the following simple result where we write $u_{\pm} = \max\{ \pm u, 0 \}$.

\begin{lemma} \label{lemma_upm}
Let $\Omega \subset \R^n$ be a bounded open set and $\alpha \in \R$. If $u \in H^1(\Omega, d^{\alpha})$, then $u_{\pm} \in H^1(\Omega, d^{\alpha})$ and the weak derivatives satisfy 
\[
\p_j u_{\pm} = \begin{cases} \p_j u & \text{ in $\{ \pm u > 0 \}$}, \\ 0 & \text{ elsewhere}. \end{cases}
\]
If $\Omega$ has Lipschitz boundary and $-1 < \alpha < 1$, we also have $T(u_{\pm}) = (Tu)_{\pm}$.
\end{lemma}
\beginproof
If $u \in H^1(\Omega, d^{\alpha})$, then it is standard that $u_{\pm} \in H^1_{\mathrm{loc}}(\Omega)$ and that $\p_j u_{\pm}$ satisfies the formula above locally in $\Omega$. It follows directly that $u_{\pm} \in H^1(\Omega, d^{\alpha})$. The formula $T(u_{\pm}) = (Tu)_{\pm}$ holds for $u \in C^{\infty}(\overline{\Omega})$, and it continues to hold for $u \in H^1(\Omega, d^{\alpha})$ by density.
\QEDA \\

\noindent {\it Proof of Theorem \ref{thm_mp}. \,}
We will prove that if $Tu \leq 0$ a.e.\ on $\p \Omega$, then $u \leq 0$ a.e.\ in $\Omega$ (the other statements follow easily from this). This will be done by testing the equation against $d^{\alpha} v$ where $v=u_+$. Let $u \in H^1(\Omega, d^{\alpha})$ and $v \in C^{\infty}_c(\Omega)$, and define the bilinear form 
\[
B(u,v) = \sum_{j,k}(\sigma^{jk} \p_j u, \p_k(d^{\alpha} v)) + (c u, d^{\alpha} v)
\]
where the inner products are in $L^2(\Omega)$. 
Using the Leibniz rule gives 
\begin{equation} \label{ajk_eq}
B(u,v) = \sum_{j,k} \lb (\sigma^{jk} d^{\alpha/2} \p_j u, d^{\alpha/2} \p_k v) + (\sigma^{jk} d^{\alpha/2} \p_j u, \alpha (\p_k d) d^{\alpha/2-1} v) + (c d^{\alpha/2} u, d^{\alpha/2} v) \rb.
\end{equation}
Now $\abs{\nabla d} \leq 1$. Using Theorem \ref{thm_trace} (b), the identity \eqref{ajk_eq} continues to hold for all $u \in H^1(\Omega, d^{\alpha})$ and $v \in H^1_0(\Omega, d^{\alpha})$ by density.

Finally, let $u$ be a solution with $Tu \leq 0$ a.e.\ on $\p \Omega$. Then $B(u,v) = 0$ for all $v \in H^1_0(\Omega, d^{\alpha})$, and $T(u_+) = 0$ by Lemma \ref{lemma_upm}. Thus we may choose $v = u_+$, which  implies that 
\[
0 = B(u, u_+) = B(u_+, u_+) - B(u_-, u_+).
\]
By Lemma \ref{lemma_upm} any product $\p^{\beta} u_+ \p^{\gamma} u_-$ vanishes a.e.\ in $\Omega$ for $\abs{\beta}, \abs{\gamma} \leq 1$. This implies that $B(u_-, u_+) = 0$, which yields $B(u_+,u_+) = 0$. We now use \eqref{ajk_eq} with $u=v=u_+$,  the assumption \eqref{ajk_elliptic} for $\sigma^{jk}$, and the assumption that $c \geq 0$ to obtain that 
\[
\lambda \norm{d^{\alpha/2} \nabla u_+}^2 \leq \Lambda \abs{\alpha} \norm{d^{\alpha/2} \nabla u_+} \, \norm{d^{\alpha/2-1} u_+}.
\]
Using the Hardy inequality from Lemma \ref{lemma_hardy}, we obtain that 
\[
\norm{d^{\alpha/2} \nabla u_+}^2 \leq \frac{\Lambda}{\lambda} \abs{\alpha} C C_1^{\alpha/2} \left(1+ \frac{1}{\abs{\alpha-1}} \right) \norm{d^{\alpha/2} \nabla u_+}^2.
\]
If $\varepsilon$ is small enough and $\abs{\alpha} \leq \varepsilon$, then the constant on the right is $\leq 1/2$. It follows that $\nabla u_+ = 0$, which implies that $u_+ = 0$ using that $Tu_+ = 0$.
\QEDA

\section{Reconstruction procedure}
\label{sec:recProc}
This section lists the reconstruction procedure for an isotropic conductivity $\sigma$ from a $2\times 2$ power density matrix $\mathbf{H}$ based on \cite{monard2012a}. One can extend this procedure for anisotropic conductivities by adding another step following the approach of \cite{BalMonard}. For simplicity, we limit ourselves to the isotropic case. Throughout this section we assume that the boundary functions $f^1$ and $f^2$ were chosen in accordance with theorem \ref{thm:Main} so that the corresponding solutions $u_1$ and $u_2$ entering $\mathbf{H}$ satisfy the non-vanishing Jacobian constraint \eqref{eq:nonvan} and are ordered so that $\text{det}[\nabla u_1 \, \nabla u_2]>0$.\\
\indent The procedure is characterized by two steps. In the first step we reconstruct the angle $\theta$ that enables us to determine the functionals $\mathbf{S}_i=\sqrt{\sigma} \nabla u_i$ from the entries of $H_{ij}=\sigma \nabla u_i \cdot \nabla u_j$. In the second step, we reconstruct $\sigma$ from the functionals $\mathbf{S}_i$.

\subsection{Reconstruction of $\theta$}
We consider the power density matrix $\mathbf{H}$ and the matrix $\mathbf{S}$ composed of the functionals $\mathbf{S}_1$ and $\mathbf{S}_2$: $\mathbf{S}=[\mathbf{S}_1 \, \mathbf{S}_2]$. By definition, $\mathbf{H}$ is symmetric and by the Jacobian constraint and the lower bound on $\sigma$ it follows that $\mathbf{H}$ is positive definite: For any $x=(x_1,x_2)\neq 0, x^T \mathbf{H} x = \sigma \abs{x_1\nabla u_1 + x_2 \nabla u_2}^2>0$ (since $\nabla u_1 $ and $\nabla u_2$ are nonzero and linearly independent by the Jacobian constraint). \\
\indent In order to split the functionals $\mathbf{S}_i$ from the entries of the power density data $H_{ij}$, $\mathbf{S}$ is orthonormalized into a $SO(2)$-valued matrix $\mathbf{R}$: $\mathbf{R}=\mathbf{S}\mathbf{T}^T$. By definition, $\mathbf{R}$ is orthogonal and has determinant one and the transfer matrix $\mathbf{T}$ is determined by the data. The question is, which matrices $\mathbf{T}$ satisfy the equality $\mathbf{R}=\mathbf{S}\mathbf{T}^T$ under the conditions on $\mathbf{R}$. This question has no unique answer, so several choices of $\mathbf{T}$ are possible, for instance $\mathbf{T}=\mathbf{H}^{-\frac{1}{2}}$ or obtaining $\mathbf{T}$ by Gram-Schmidt orthonormalization. As $\mathbf{R}$ is a rotation matrix, it is parameterized by the angle function $\theta$ as follows:
\begin{equation*}
    \mathbf{R}(\theta) = \begin{bmatrix} \cos \theta & -\sin \theta\\ \sin \theta & \cos \theta \end{bmatrix}.
\end{equation*}
From this definition, we see that once $\mathbf{T}$ and $\mathbf{S}$ are known, the function $\theta$ can be computed by
\begin{equation*}
    \theta=\text{arg}(\mathbf{R}_1),
\end{equation*}
where $\mathbf{R}_1$ denotes the first column of $\mathbf{R}$. The orthonormalization technique and thus the choice of $\mathbf{T}$ influences the angle $\theta$. Our choice of $\mathbf{T}$ and the corresponding interpretation of $\theta$ is discussed in subsection \ref{sec:ChoiceTtheta}. \par
Defining $\mathbf{T}=(T_{ij})_{1\leq i,j \leq 2}$ and $\mathbf{T}^{-1}=(T^{ij})_{1\leq i,j \leq 2}$, and letting 
\begin{equation*}
    \mathbf{V_{ij}}=\nabla (T_{i1})T^{1j}+\nabla (T_{i2})T^{2j},
\end{equation*}
then $\theta$ is determined by the following equation \cite[Eq. (65)]{monard2012a}:
\begin{equation}\label{eq:theta}
    \nabla \theta = \mathbf{F},
\end{equation}
with
\begin{equation*}
    \mathbf{F}=\frac{1}{2}(\mathbf{V_{12}}-\mathbf{V_{21}}-\mathcal{J}\nabla \log D),
\end{equation*}
$\mathcal{J}=\begin{bmatrix} 0 & -1\\ 1 & 0\end{bmatrix}$, and $D=(H_{11}H_{22}-H_{12}^2)^{\frac{1}{2}}$. Once $\theta$ is known at at least one point on the boundary one can integrate $\mathbf{F}$ along curves originating from that point to obtain $\theta$ throughout the whole domain. Alternatively, when assuming that $\theta$ is known along the whole boundary one can apply the divergence operator to \eqref{eq:theta} and solve the following Poisson equation with Dirichlet boundary condition:
\begin{equation}\label{eq:thetaBVP}
    \begin{cases}
        \Delta \theta=\nabla \cdot \mathbf{F} & \text{in }\Omega,\\
        \theta=\theta_{\text{true}} & \text{on }\partial \Omega.
    \end{cases}
\end{equation}
In our implementation, we use the second option and discuss in subsection \ref{sec:ChoiceTtheta} knowledge of $\theta$ along the boundary.

\subsection{Reconstruction of $\sigma$}
Reconstruction of $\sigma$ is based on \cite[Eq. (68)]{monard2012a}
\begin{equation}\label{eq:sigmaRec}
    \nabla \log \sigma = \mathbf{G},
\end{equation}
with 
\begin{align*}
    \mathbf{G}&=\cos(2 \theta)\mathbf{K} + \sin(2\theta)\mathbf{K},\\
    \mathbf{K}&=\mathcal{U}(\mathbf{V_{11}}-\mathbf{V_{22}})+\mathcal{J}\mathcal{U}(\mathbf{V_{12}}-\mathbf{V_{21}}) \quad \text{and} \quad \mathcal{U}=\begin{bmatrix} 1 & 0\\ 0 & -1\end{bmatrix}.
\end{align*}
Similarly as for $\theta$ one need to solve a gradient equation to obtain $\sigma$ and has the possibility of either integrating along curves or solving a Poisson equation, assuming knowledge of $\sigma$ in one point or along the whole boundary respectively. We assume knowledge of $\sigma$ along the whole boundary and solve the following Poisson equation with Dirichlet condition:
\begin{equation}\label{eq:sigmaBVP}
    \begin{cases}
        \Delta \log (\sigma)= \nabla \cdot \mathbf{G} &\text{in }\Omega,\\
        \log (\sigma)=\log (\sigma_{\text{true}}) & \text{on }\partial \Omega.
    \end{cases}
\end{equation}

\subsection{Choice of the transfer matrix $\mathbf{T}$ and knowledge of $\theta$}
\label{sec:ChoiceTtheta}

For our implementation, we use Gram-Schmidt orthonormalization to obtain the transfer matrix $\mathbf{T}$:
\begin{equation*}
    \mathbf{T}=\begin{bmatrix} H_{11}^{-\frac{1}{2}} & 0 \\ -H_{12}H_{11}^{-\frac{1}{2}}D^{-1} & H_{11}^{\frac{1}{2}}D^{-1}\end{bmatrix}.
\end{equation*}
By the Jacobian constraint, it follows that $H_{11}>0$ so that $\mathbf{T}$ is well-defined. As a direct consequence of using Gram-Schmidt orthonormalization the first column of $\mathbf{R}$ simplifies to:
\begin{equation*}
    \mathbf{R}_1 = T_{11}\mathbf{S}_1 + T_{12} \mathbf{S}_2 = \frac{\cancel{\sqrt{\sigma}}\nabla u_1}{\cancel{\sqrt{\sigma}}\abs{\nabla u_1}}.
\end{equation*}
Therefore, the angle $\theta$ simply defines the angle between $\nabla u_1$ and the $x_1$-axis. Hence,
\begin{equation}\label{eq:thetau1}
    \theta = \text{arg}(\nabla u_1).
\end{equation}
In addition, following this definition for $\mathbf{T}$ the vector fields $\mathbf{V_{ij}}$ can be written explicitly in terms of $\mathbf{H}$:
\begin{align}\label{eq:Vexp}
\begin{aligned}
    \mathbf{V_{11}}&= \nabla \log H_{11}^{-\frac{1}{2}}, & \mathbf{V_{12}}&=0,\\ \mathbf{V_{21}}&= -\frac{H_{11}}{D} \nabla \left( \frac{H_{12}}{H_{11}}\right), & \mathbf{V_{22}}&=\nabla \log \left( \frac{H_{11}^{\frac{1}{2}}}{D}\right).
    \end{aligned}
\end{align}
Knowledge of $\theta$ at the boundary is essential for the reconstruction procedure. By this definition of $\mathbf{T}$, knowledge of $\theta$ is directly related to knowledge of the gradient $\nabla u_1$ and the current $\sigma \nabla u_1$, as both vector fields have the same direction. We decompose $\sigma \nabla u_1$ into two parts with contribution from the unit outward normal $\nu$ and the tangent vector $\omega=\mathcal{J}\nu$:
\begin{equation*}
    \sigma \nabla u_1 = (\sigma \nabla u_1 \cdot \nu)\nu + (\sigma \nabla u_1 \cdot \omega)\omega.
\end{equation*}
As along the whole boundary $u_1$ and $\sigma$ are known, the contribution from $\sigma \nabla u_1 \cdot \omega$ is known as well. Furthermore, along the part of the boundary $\partial \Omega \backslash \Gamma$ we have additional information as $u_1$ vanishes. Therefore, the only contribution is from the unit outward normal $\nu$, so that on this part of the boundary $\sigma \nabla u_1$ has either the direction of $\nu$, $-\nu$ or the zero vector. However, in order to have full information of $\theta$ along the boundary one needs knowledge about the Neumann data $\sigma \nabla u_1 \cdot \nu$.

\section{Numerical Examples}
\begin{sloppypar}
The \textsc{Matlab} and \textsc{Python} code to generate the numerical examples can be found on \textsc{GitLab}: \href{https://lab.compute.dtu.dk/hjsc/jacobian-of-solutions-to-the-conductivity-equation-in-limited-view.git}{https://lab.compute.dtu.dk/hjsc/jacobian-of-solutions-to-the-conductivity-equation-in-limited-view.git}.\par
\begin{figure}[ht!]
    \centering
    \begin{minipage}[t]{0.5\textwidth}
        \centering
        \includegraphics[width=\linewidth]{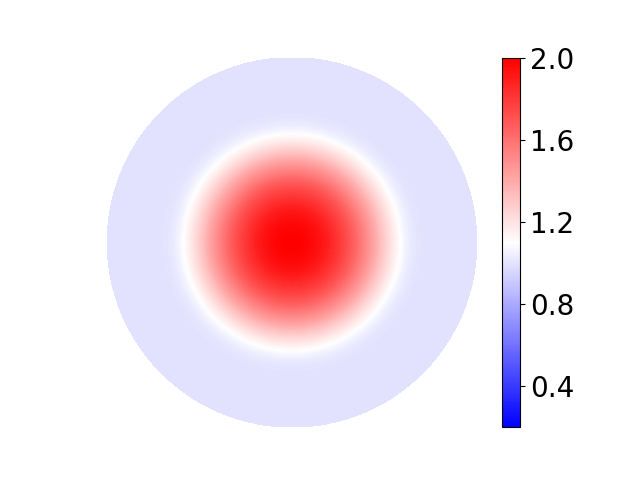}
        \caption*{Case 1}
    \end{minipage}%
    \begin{minipage}[t]{0.5\textwidth}
        \centering
        \includegraphics[width=\linewidth]{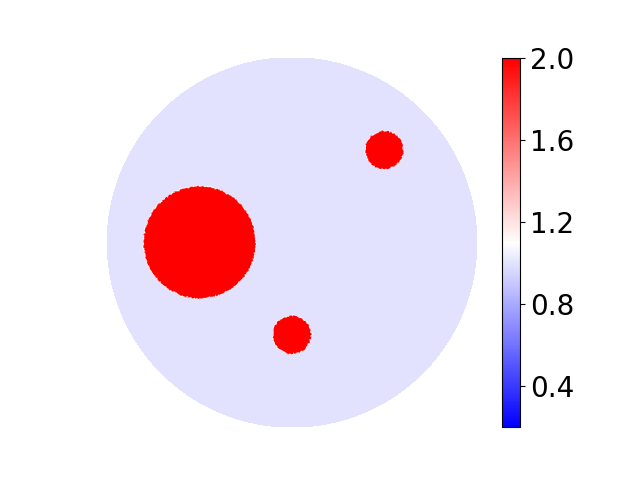}
        \caption*{Case 2}
    \end{minipage}
    \caption{The conductivities $\sigma$ used for the reconstruction procedure.}
    \label{fig:TrueSig}
\end{figure}
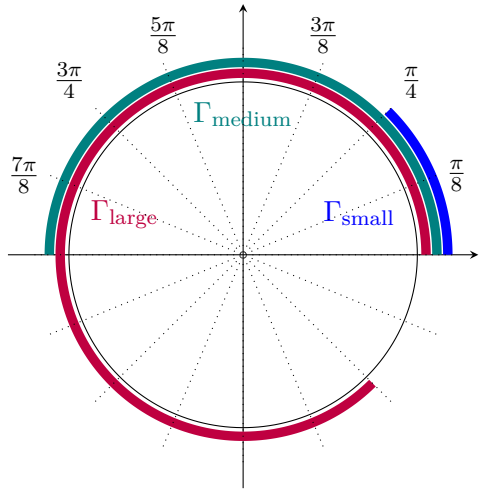
\begin{figure}[ht!]
    \centering
    \begin{minipage}[t]{0.5\textwidth}
        \centering
        \begin{tikzpicture}
\begin{axis}[
	trig format plots=rad,
	axis x line=center,
    axis y line=center,
    xlabel={},
    ylabel={},
    ticks=none,
    unit vector ratio*=1 1 1,
    xmin=-1.35, xmax=1.35,
    ymin=-1.35, ymax=1.45,
    height=8cm,
]
\tikzmath{\v=pi/8;} 

\addplot[domain=0:2*pi, samples=200] ({cos(x)},{sin(x)});

\addplot[domain=0:7*pi/4, samples=200, purple, line width=.75ex] ({1.05*cos(x)},{1.05*sin(x)}) node [pos=.5,xshift=4ex,yshift=0ex,anchor=north] {$ \Gamma_{\textup{large}}$};
\addplot[domain=0:pi, samples=200, teal, line width=.75ex] ({1.1125*cos(x)},{1.1125*sin(x)}) node [pos=.5,xshift=0ex, yshift=-2ex,anchor=north] {$ \Gamma_{\textup{medium}}$};
\addplot[domain=0:pi/4, samples=200, blue, line width=.75ex] ({1.175*cos(x)},{1.175*sin(x)}) node [pos=.5,xshift=-1.8ex,yshift=-5ex,anchor=south east] {$ \Gamma_{\textup{small}}$};

\addplot[domain=-1.2:1.2, dotted] ({cos(\v)*x},{sin(\v)*x}) node [anchor=west] {$ \frac{\pi}{8}$};
\addplot[domain=-1.2:1.2, dotted] ({cos(2*\v)*x},{sin(2*\v)*x}) node [anchor=south west] {$ \frac{\pi}{4}$};
\addplot[domain=-1.2:1.2, dotted] ({cos(3*\v)*x},{sin(3*\v)*x}) node [anchor=south] {$ \frac{3\pi}{8}$};
\addplot[domain=-1.2:1.2, dotted] ({cos(4*\v)*x},{sin(4*\v)*x});
\addplot[domain=-1.2:1.2, dotted] ({cos(5*\v)*x},{sin(5*\v)*x}) node [anchor=south] {$ \frac{5\pi}{8}$};
\addplot[domain=-1.2:1.2, dotted] ({cos(6*\v)*x},{sin(6*\v)*x}) node [anchor=south east] {$ \frac{3\pi}{4}$};
\addplot[domain=-1.2:1.2, dotted] ({cos(7*\v)*x},{sin(7*\v)*x}) node [anchor=east] {$ \frac{7\pi}{8}$};
\addplot[domain=-1.2:1.2, dotted] ({cos(8*\v)*x},{sin(8*\v)*x});

\end{axis}
\end{tikzpicture}
    \end{minipage}
    \caption{Different sizes of $\Gamma$ used for the reconstruction procedure.}
    \label{fig:GammaSize}
\end{figure}

\end{sloppypar}
Our aim is to illustrate numerically how two boundary conditions can be selected so that the non-vanishing Jacobian condition \eqref{eq:nonvan} for corresponding solutions is satisfied in accordance with theorem \ref{thm:Main}. And we choose the order of the corresponding solutions so that $\mathrm{det}[\nabla u_1 \, \nabla u_2]>0$. Furthermore, we show numerically how this can be used to reconstruct the conductivity from power density data. For that purpose, we implemented the reconstruction procedure in section \ref{sec:recProc} in \python{} and used \fenics{}~\cite{fenics} to solve the PDEs. We use a fine mesh to generate our power density data and a coarser mesh to address the reconstruction problem. We use $N_{\text{data}}=79281$ nodes in the high-resolution case, while for the coarser mesh we consider a resolution of $N_{\text{recon}}=50845$ nodes. For both meshes, we use $\mathbb{P}_1$ elements. We consider the domain $\Omega$ to be the unit disk: $\Omega = B(\mathbf{0},1)$. Furthermore, we consider two test cases for an isotropic conductivity $\sigma$ defined by:
\begin{align*}
    \sigma_{\text{case 1}}(x_1,x_2) &= \begin{cases}1 + e^{\left(2\,-\,\frac{2}{1-\frac{(x_1)^2+(x_2)^2}{1-0.8^2}}\right)} & 0\leq (x_1)^2+(x_2)^2 \leq 0.8^2,\\
    1 & 0.8^2 \leq (x_1)^2+(x_2)^2 \leq 1,
    \end{cases}\\
    \sigma_{\text{case 2}}(x_1,x_2) &= \begin{cases}2 & \left(x_1+\frac{1}{2} \right)^2 + (x_2)^2 \leq 0.3^2,\\
        2 & \left(x_1\right)^2 + \left(x_2+\frac{1}{2} \right)^2 \leq 0.1^2,\\
        2 & \left(x_1-\frac{1}{2} \right)^2 + \left(x_2-\frac{1}{2} \right)^2 \leq 0.1^2,\\
        1 & \text{otherwise}.
    \end{cases}
\end{align*}
for $(x_1,x_2)\in \Omega$. Figure \ref{fig:TrueSig} illustrates the conductivities. To investigate influence of the size of the boundary of control, $\Gamma$, we consider three different sizes that are outlined in figure \ref{fig:GammaSize}. \par

\begin{figure}[ht!]
    \centering
    \includegraphics[width=0.8\textwidth]{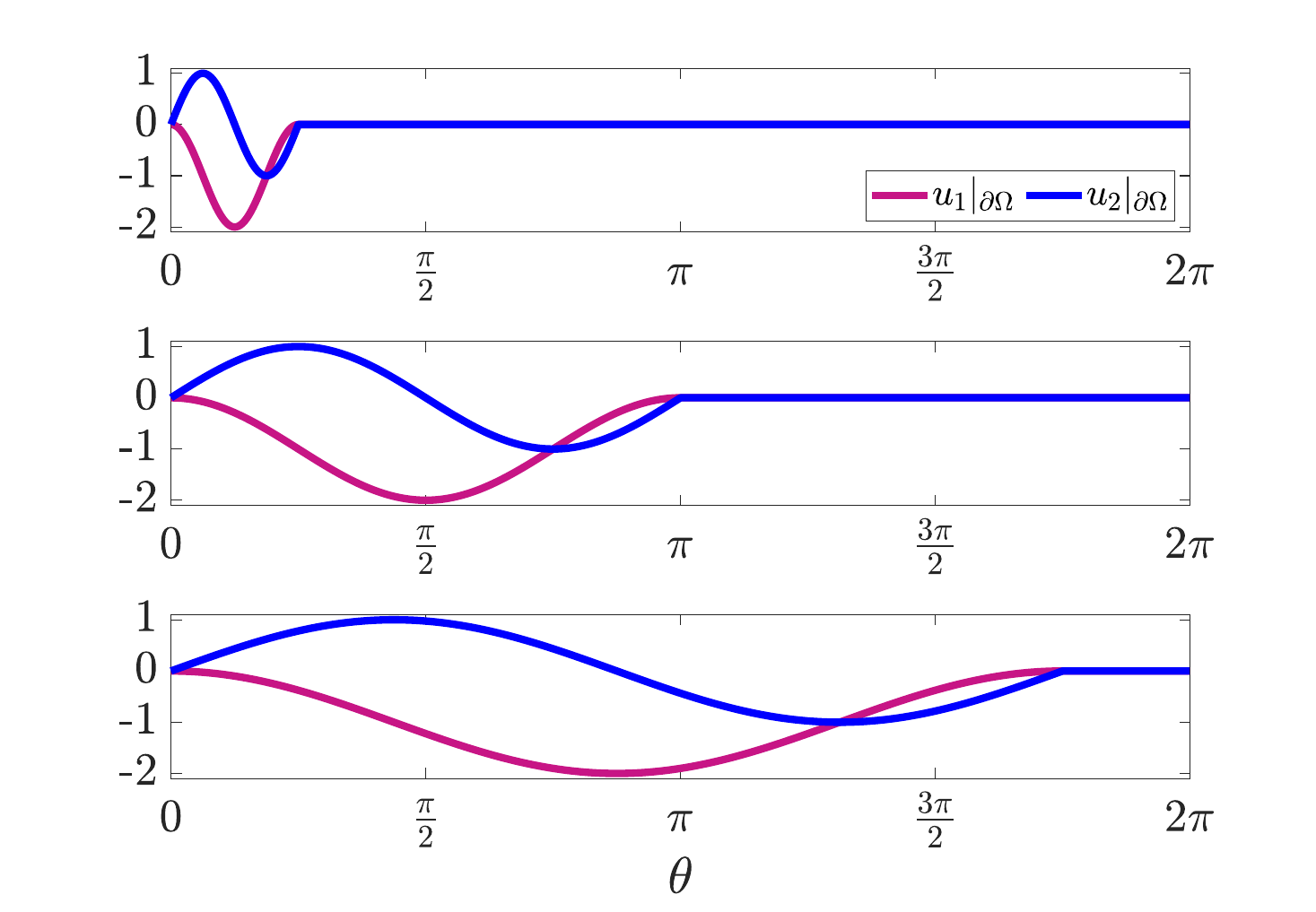}
    \caption{The continuous boundary functions $u_1^c\vert_{\partial \Omega}$ and $u_2^c\vert_{\partial \Omega}$ used for the reconstruction procedure for $\Gamma_{\text{small}}$, $\Gamma_{\text{medium}}$ and $\Gamma_{\text{large}}$ (top to bottom).}
    \label{fig:contBC}
\end{figure}

We demonstrate that the Jacobian constraint is satisfied for a choice of continuous and discontinuous boundary conditions in accordance with theorem \ref{thm:Main}.
The functions $(f_1^c,f_2^c)$ below yield continuous boundary functions that satisfy condition (a) in theorem \ref{thm:Main}: From the right part in figure \ref{fig:ws} we see that $\mathrm{arg}(\dot{\gamma}^c)$ is strictly increasing and thus monotone and $\dot{\gamma^c}$ satisfies $\abs{\mathrm{Ind}(\dot{\gamma}^c)}=1$. As $\mathrm{Ind}(\dot{\gamma}^c)$ denotes the winding number of $\dot{\gamma}^c$ it can also be observed visually from the left part of figure \ref{fig:ws} that its winding number is 1.

\begin{equation*}
   (f_1^c(\eta(t)),f_2^c(\eta(t))) = \begin{cases}
        (\cos(8t)-1,\sin(8t)) & \text{for }\Gamma_{\text{small}}=\{t \in \left[0,\frac{\pi}{4}\right]\}\\
        (\cos(2t)-1,\sin(2t)) & \text{for }\Gamma_{\text{medium}}=\{t \in \left[0,\pi\right]\}\\
        \left(\cos \left(\frac{8t}{7}\right)-1,\sin \left(\frac{8t}{7}\right)\right) & \text{for }\Gamma_{\text{large}}=\left\{t \in \left[0,\frac{7\pi}{4}\right] \right\}.
    \end{cases}
\end{equation*}
The corresponding functions $u_i^c\vert_{\partial \Omega}$ extended by zero along the whole boundary are illustrated in figure \ref{fig:contBC}.

The functions $(f_1^d,f_2^d)$ below yield discontinuous boundary functions that satisfy condition (b) in theorem \ref{thm:Main}: From the right part in figure \ref{fig:ws} we see that $\mathrm{arg}(\dot{\gamma}^d)$ is strictly increasing and thus monotone and $\dot{\gamma}^d$ satisfies $\abs{\mathrm{Ind}(\dot{\gamma}^d)}=1/2$. As we interpret $\mathrm{Ind}(\dot{\gamma}^d)$ as the winding number of $\dot{\gamma}^d$ it can also be observed visually from the left part of figure \ref{fig:ws} that its winding number is 1/2.

\begin{equation*}
   (f_1^d(\eta(t)),f_2^d(\eta(t))) = \begin{cases}
        (\cos(4t)-1,\sin(5t)) & \text{for }\Gamma_{\text{small}}=\{t \in \left[0,\frac{\pi}{4}\right]\}\\
        \left(\cos \left(t\right)-1,\sin\left(\frac{5t}{4}\right)\right) & \text{for }\Gamma_{\text{medium}}=\{t \in \left[0,\pi\right]\}\\
        \left(\cos \left(\frac{4t}{7}\right)-1,\sin\left(\frac{5t}{7}\right)\right) & \text{for }\Gamma_{\text{large}}=\left\{t \in \left[0,\frac{7\pi}{4}\right] \right\}.
    \end{cases}
\end{equation*}

\begin{figure}[ht!]
    \centering
    \includegraphics[width=0.8\textwidth]{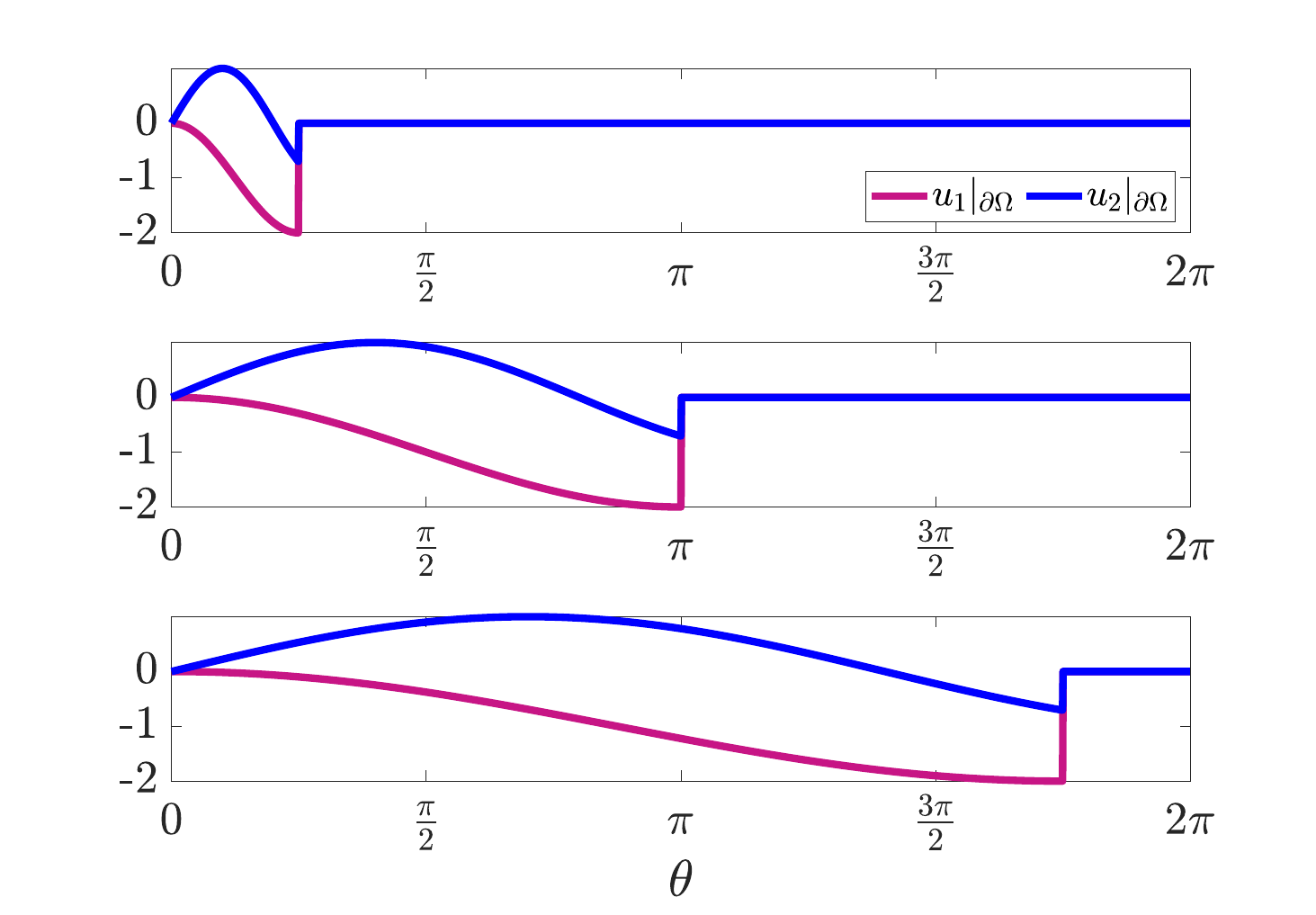}
    \caption{The discontinuous boundary functions $u_1^d\vert_{\partial \Omega}$ and $u_2^d\vert_{\partial \Omega}$ used for the reconstruction procedure for $\Gamma_{\text{small}}$, $\Gamma_{\text{medium}}$ and $\Gamma_{\text{large}}$ (top to bottom).}
    \label{fig:discontBC}
\end{figure}

\begin{figure}[ht!]
    \centering
    \begin{minipage}[t]{0.5\textwidth}
        \centering
        \includegraphics[width=\linewidth]{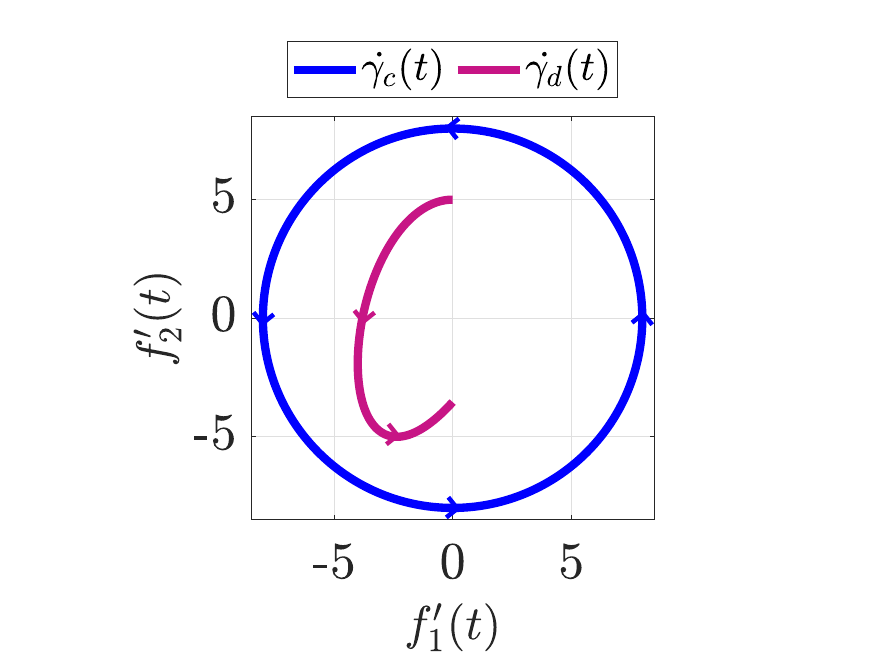}
    \end{minipage}%
    \begin{minipage}[t]{0.5\textwidth}
        \centering
        \includegraphics[width=\linewidth]{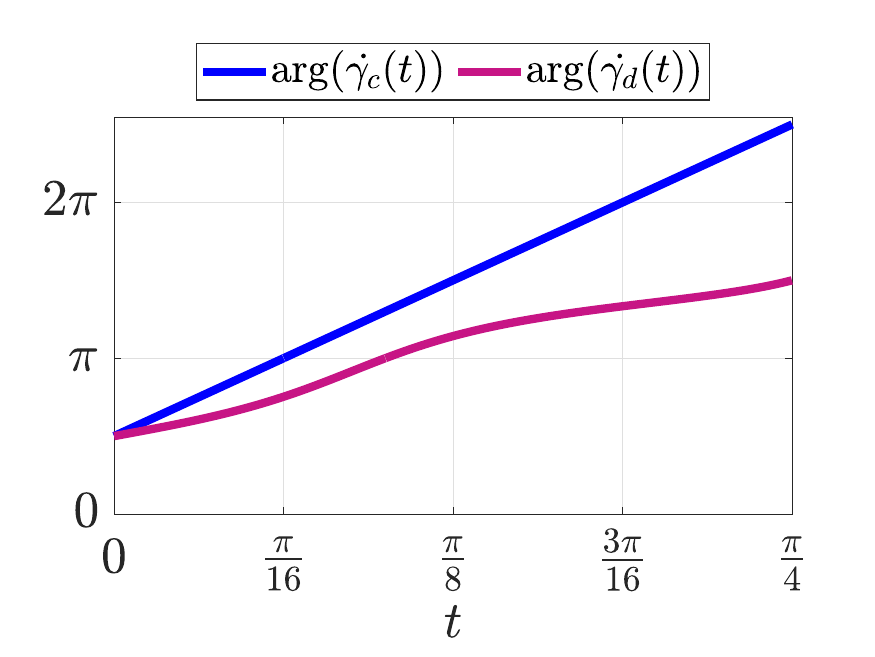}
    \end{minipage}
     \caption{ Illustration that the functions in \ref{fig:contBC} and \ref{fig:discontBC} respectively satisfy the conditions in theorem \ref{thm:Main}. The left part illustrates the winding number $\mathrm{Ind}(\dot{\gamma})$ of the curves $\dot{\gamma_c}(t)=((f_1^c)'(t),(f_2^c)'(t))$ and $\dot{\gamma_d}(t)=((f_1^d)'(t),((f_2^d)'(t))$, while the right part shows the behavior of their arguments.} 
    \label{fig:ws}
\end{figure}

The corresponding functions $u_i^d\vert_{\partial \Omega}$ extended by zero along the whole boundary are illustrated in figure \ref{fig:discontBC}.

\begin{figure}[ht!]
    \centering
    \begin{minipage}[t]{0.03\textwidth}
        \begin{turn}{90}\hspace{3mm}Continuous $u_1\vert_{\partial  \Omega}$\end{turn}
    \end{minipage}%
    \begin{minipage}[t]{0.32\textwidth}
        \centering
        \includegraphics[trim={0cm 0cm 4cm 0cm},clip,width=0.76\linewidth]{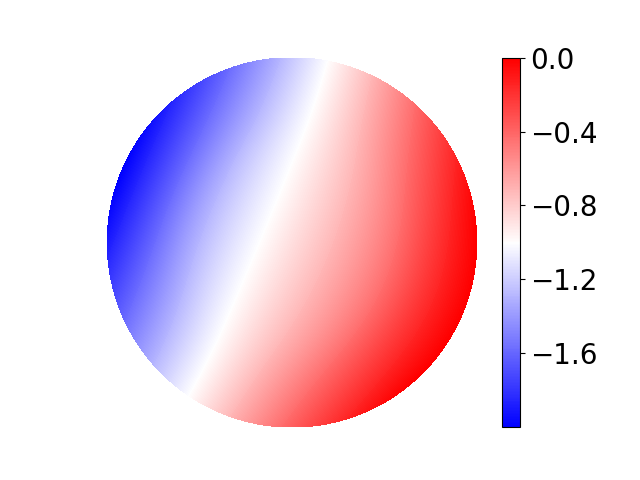}
    \end{minipage}%
    \begin{minipage}[t]{0.32\textwidth}
        \centering
        \includegraphics[trim={0cm 0cm 4cm 0cm},clip,width=0.76\linewidth]{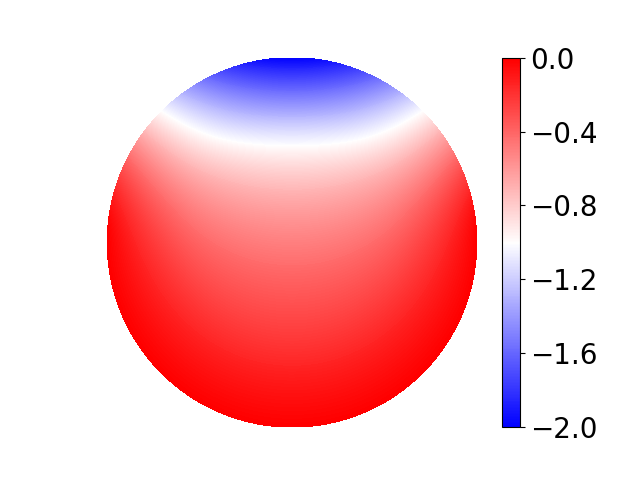}
    \end{minipage}%
    \begin{minipage}[t]{0.32\textwidth}
        \centering
        \includegraphics[width=\linewidth]{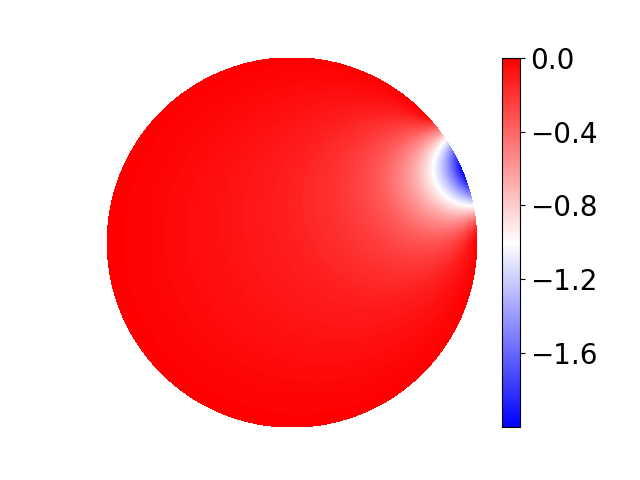}
    \end{minipage}
    \begin{minipage}[t]{0.03\textwidth}
        \begin{turn}{90}\hspace{3mm}Discontinuous $u_1\vert_{\partial  \Omega}$\end{turn}
    \end{minipage}%
    \begin{minipage}[t]{0.32\textwidth}
        \centering
        \includegraphics[trim={0cm 0cm 4cm 0cm},clip,width=0.76\linewidth]{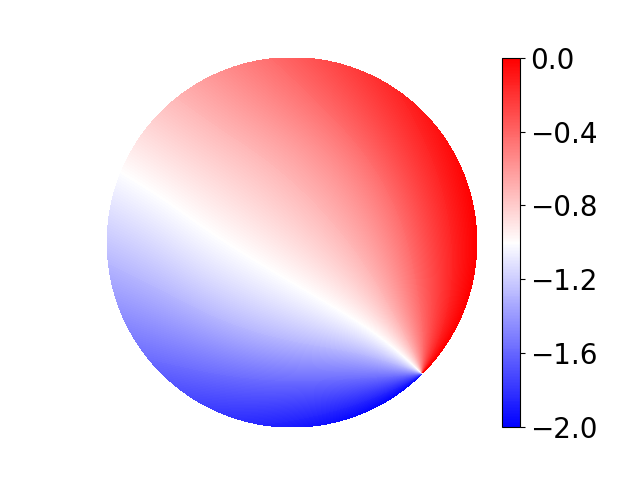}
    \end{minipage}%
    \begin{minipage}[t]{0.32\textwidth}
        \centering
        \includegraphics[trim={0cm 0cm 4cm 0cm},clip,width=0.76\linewidth]{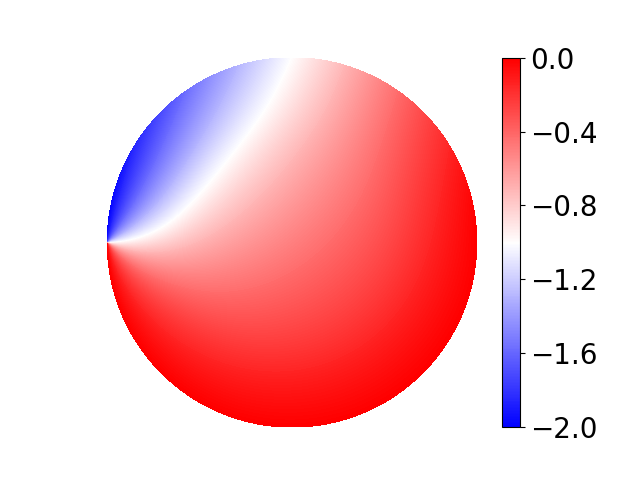}
    \end{minipage}%
    \begin{minipage}[t]{0.32\textwidth}
        \centering
        \includegraphics[width=\linewidth]{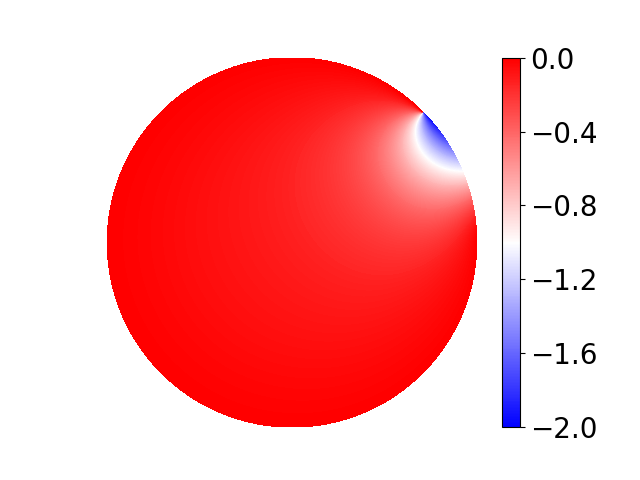}
    \end{minipage}
    \begin{minipage}[t]{0.03\textwidth}
        \begin{turn}{90}\hspace{3mm}Continuous $u_2\vert_{\partial  \Omega}$\end{turn}
    \end{minipage}%
    \begin{minipage}[t]{0.32\textwidth}
        \centering
        \includegraphics[trim={0cm 0cm 4cm 0cm},clip,width=0.76\linewidth]{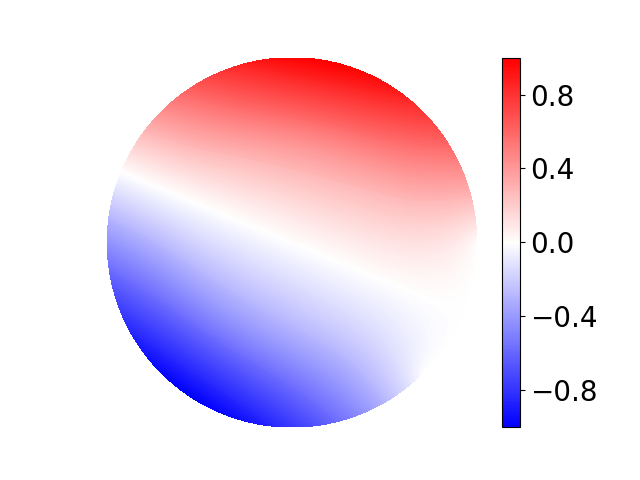}
    \end{minipage}%
    \begin{minipage}[t]{0.32\textwidth}
        \centering
        \includegraphics[trim={0cm 0cm 4cm 0cm},clip,width=0.76\linewidth]{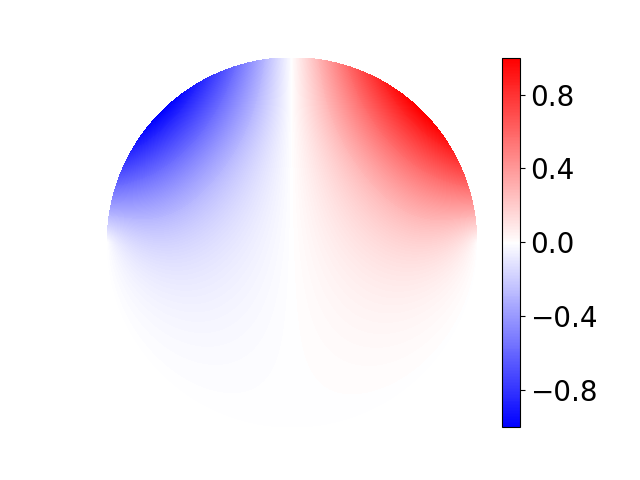}
    \end{minipage}%
    \begin{minipage}[t]{0.32\textwidth}
        \centering
        \includegraphics[width=\linewidth]{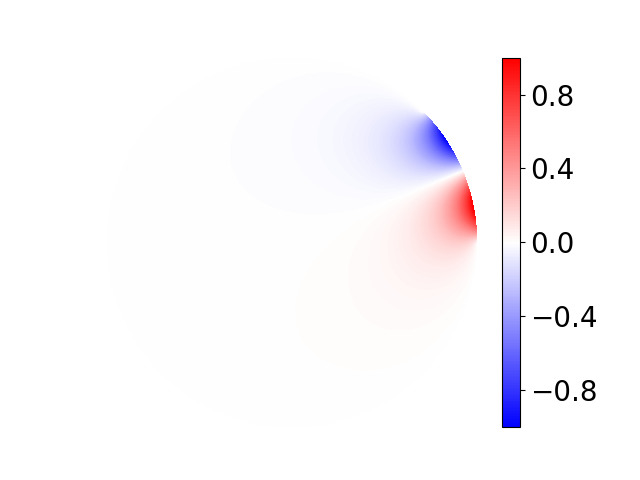}
    \end{minipage}
    \begin{minipage}[t]{0.03\textwidth}
        \begin{turn}{90}\hspace{3mm}Discontinuous $u_2\vert_{\partial  \Omega}$\end{turn}
    \end{minipage}%
    \begin{minipage}[t]{0.32\textwidth}
        \centering
        \includegraphics[trim={0cm 0cm 4cm 0cm},clip,width=0.76\linewidth]{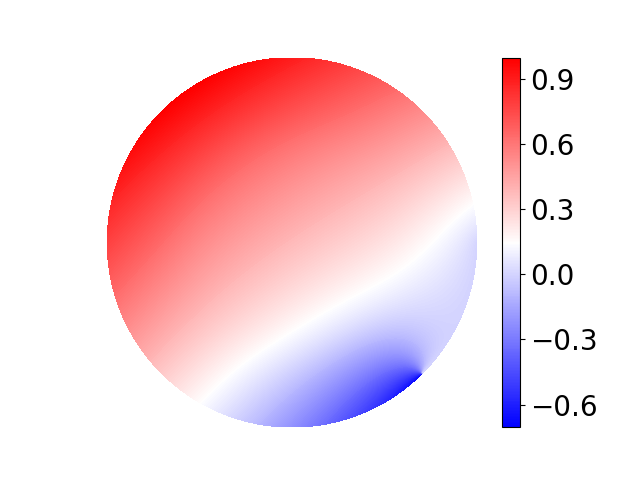}
        \caption*{\hspace{1cm}$\Gamma_{\text{large}}$}
    \end{minipage}%
    \begin{minipage}[t]{0.32\textwidth}
        \centering
        \includegraphics[trim={0cm 0cm 4cm 0cm},clip,width=0.76\linewidth]{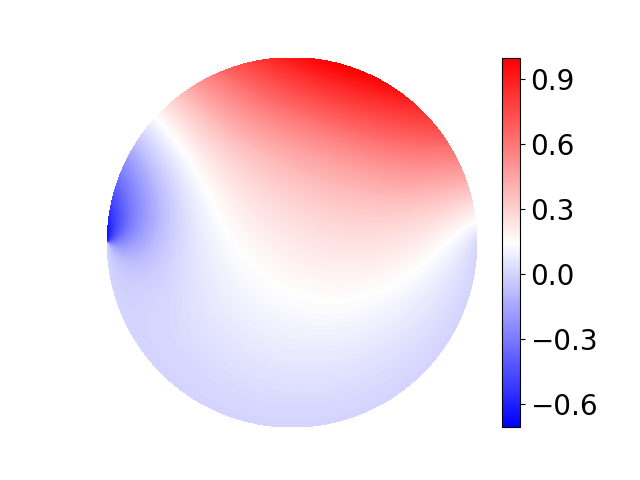}
        \caption*{\hspace{1cm}$\Gamma_{\text{medium}}$}
    \end{minipage}%
    \begin{minipage}[t]{0.32\textwidth}
        \centering
        \includegraphics[width=\linewidth]{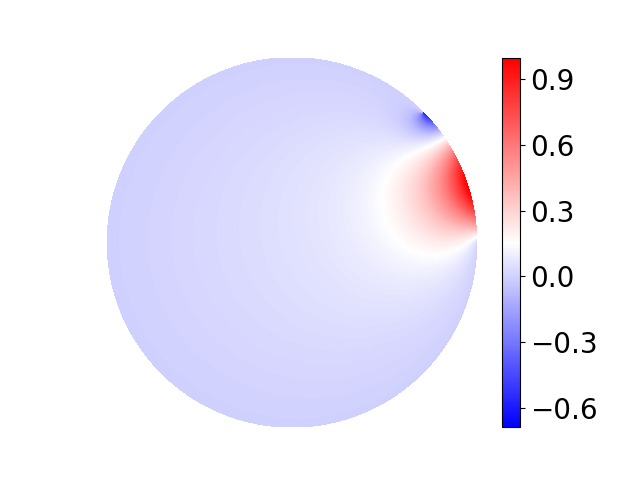}
        \caption*{$\Gamma_{\text{small}}$}
    \end{minipage}
    \caption{The solutions $u_1$ and $u_2$ induced by the discontinuous and continuous boundary conditions for $\sigma$ as in case 1 and varying boundaries of control $\Gamma_{\text{large}}$, $\Gamma_{\text{medium}}$ and $\Gamma_{\text{small}}$.}
    \label{fig:us}
\end{figure}

\begin{figure}[ht!]
    \centering
    \begin{minipage}[t]{0.03\textwidth}
        \begin{turn}{90}\hspace{6mm}Continuous BC\end{turn}
    \end{minipage}%
    \begin{minipage}[t]{0.32\textwidth}
        \centering
        \includegraphics[width=\linewidth]{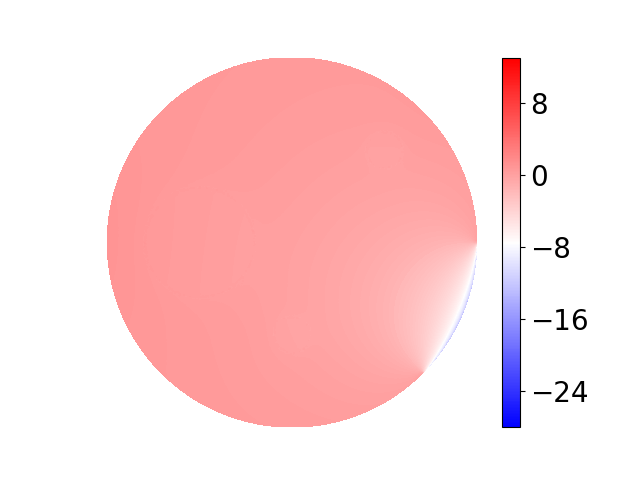}
    \end{minipage}%
    \begin{minipage}[t]{0.32\textwidth}
        \centering
        \includegraphics[width=\linewidth]{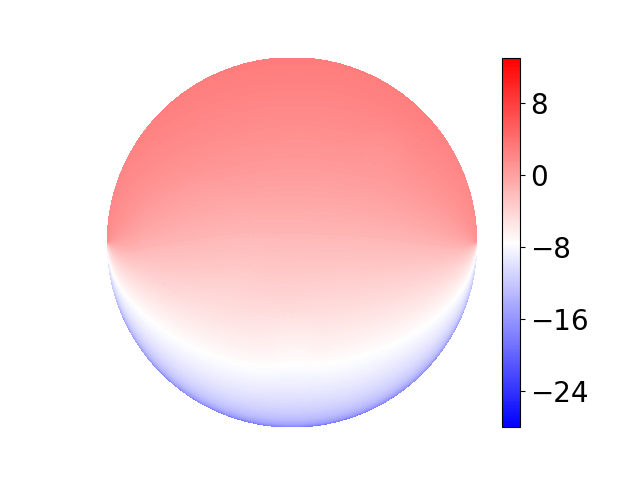}
    \end{minipage}%
    \begin{minipage}[t]{0.32\textwidth}
        \centering
        \includegraphics[width=\linewidth]{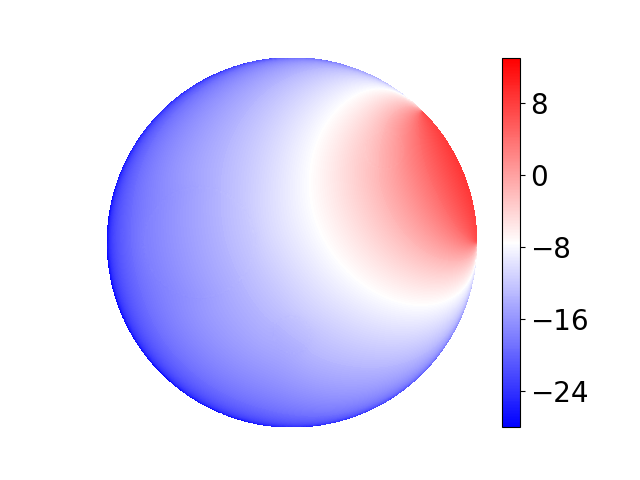}
    \end{minipage}
    \begin{minipage}[t]{0.03\textwidth}
        \begin{turn}{90}\hspace{3mm}Discontinuous BC\end{turn}
    \end{minipage}%
    \begin{minipage}[t]{0.32\textwidth}
        \centering
        \includegraphics[width=\linewidth]{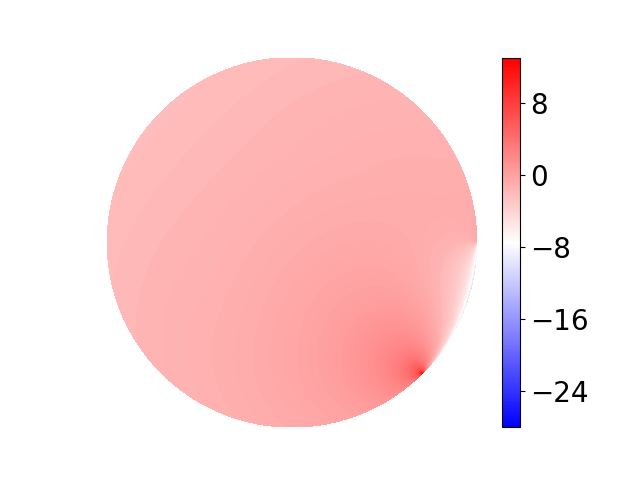}
        \caption*{$\Gamma_{\text{large}}$}
    \end{minipage}%
    \begin{minipage}[t]{0.32\textwidth}
        \centering
        \includegraphics[width=\linewidth]{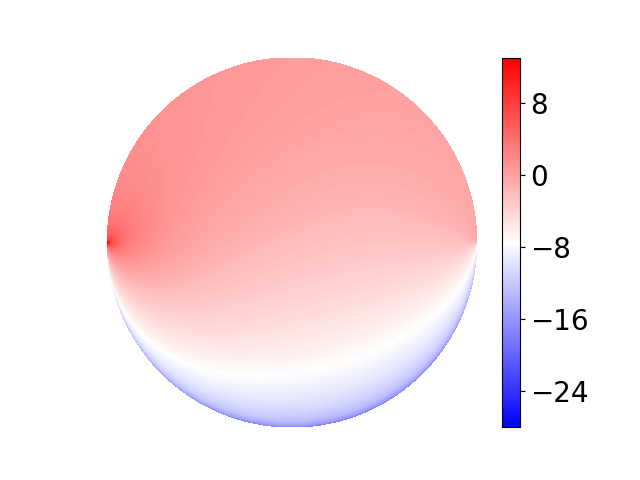}
        \caption*{$\Gamma_{\text{medium}}$}
    \end{minipage}%
    \begin{minipage}[t]{0.32\textwidth}
        \centering
        \includegraphics[width=\linewidth]{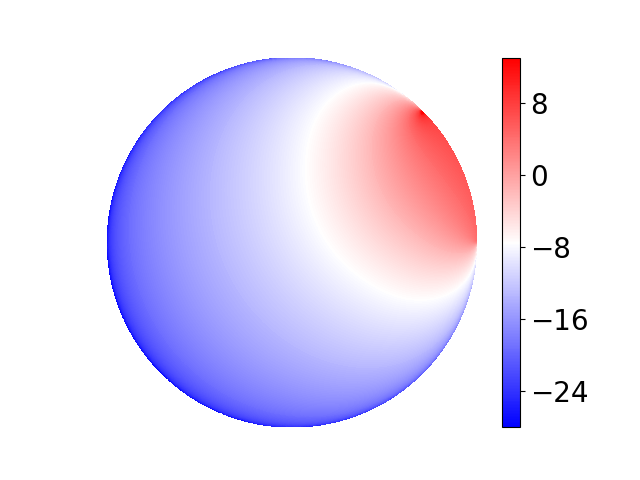}
        \caption*{$\Gamma_{\text{small}}$}
    \end{minipage}
     \caption{The expression $\log(\mathrm{det}(\mathbf{H}))$ for varying sizes of $\Gamma$, when using continuous and discontinuous boundary conditions. Large negative values (blue regions) correspond to values of $\mathrm{det}(\mathbf{H})$ close to zero.}
    \label{fig:Jac}
\end{figure}

We compute the corresponding solutions $u_1$ and $u_2$ and the three power densities $H_{11}, H_{12}$ and $H_{22}$. It is not straightforward to compute the solutions numerically for discontinuous boundary conditions; therefore, the procedure is discussed in section \ref{sec:DBC}. The solutions $u_1$ and $u_2$ are then illustrated in figure \ref{fig:us}. From table \ref{tab:RelErrSize} we see that the Jacobian condition is satisfied for all cases, as the determinant of $\mathbf{H}$ is positive. However, for a small boundary of control the values of the determinant are very small, so for $\Gamma_{\text{small}}$ the minimum values are of order $10^{-14}$. To investigate where the small values of $\mathrm{det}(\mathbf{H})$ concentrate, we illustrate the expression $\log(\mathrm{det}(\mathbf{H}))$ for continuous and discontinuous boundary conditions in figure \ref{fig:Jac}. The expression $\log(\mathrm{det}(\mathbf{H}))$ will give us more nuances of the small values than the determinant itself. From the figure, it is evident that the small values concentrate close to the boundary that we cannot control. Furthermore, the figure shows that for discontinuous boundary conditions there appear larger values of $\log(\mathrm{det}(\mathbf{H}))$ than for the continuous case, but these are mainly concentrated around the discontinuity. For the continuous boundary conditions, the maximal values are smaller, but they are more evenly distributed along the boundary of control. As the Jacobian condition is satisfied, we can use the reconstruction procedure outlined in section \ref{sec:recProc} to reconstruct the two conductivities $\sigma_{\text{case 1}}$ and $\sigma_{\text{case 2}}$. For the reconstruction procedure, we use knowledge of the true angle $\theta$ that can be computed by knowledge of the true gradient $\nabla u_1$.\\

\subsection{Solving the conductivity equation numerically with discontinuous boundary functions $u_i^d\vert_{\partial \Omega}$}
\label{sec:DBC}
It is a challenge numerically to use \fenics{} to compute the solutions with discontinuous boundary conditions as the solutions are only in $H^1(\Omega, d^{1-2s})$ and not in $H^1(\Omega)$. Using Lagrange basis functions one can define a $H^1$-function space, but this does not allow for discontinuities. To allow for discontinuities one has the possibility to define a function space using discontinuous Galerkin basis functions between the nodes, but in this way one loses interior regularity as well. Both possibilities are not optimal as in our case we only have a discontinuity at the boundary, while away from the boundary the function behaves like an $H^1$-function. For that purpose, we consider the functions $w_i=u_i-u_i^0$, were $u_i^0 \in H^1(\Omega, d^{1-2s})$ solves the Laplace equation with boundary conditions:
\begin{equation}\label{eq:LapBVP}
    \begin{cases}
        \Delta u_i^0 = 0 & \text{in }\Omega,\\
        u_i^0 = f_i & \text{on }\Gamma\\
        u_i^0 = 0 & \text{on }\partial \Omega \backslash \Gamma.
    \end{cases}
\end{equation}
Now $w^i$ solves the following boundary value problem:
\begin{equation}\label{eq:weq}
\begin{cases}
    -\text{div}(\sigma \nabla w_i) = \text{div}(\sigma \nabla u_i^0) & \text{in }\Omega,\\
    w_i = 0 & \text{on }\partial \Omega.
    \end{cases}
\end{equation}
As we consider conductivities $\sigma$ that are one on and in a neighborhood of the boundary, the right hand side $\text{div}(\sigma \nabla u_i^0)$ vanishes in a neighborhood where the discontinuity appears. Using these choices of $\sigma$ thus ensure that the discontinuity is covered so that the right hand side satisfies $\text{div}(\sigma \nabla u_i^0) \in H^{-1}(\Omega)$ implying that $w_i$ is a solution in $H^1(\Omega)$. We solve the boundary value problem \eqref{eq:LapBVP} semi-analytically for $u_i^0$ in \matlab{} using the Fourier transform. This gives us the exact solution at each node apart from the Gibbs phenomenon happening at the discontinuity. Afterwards we solve the boundary value problem \eqref{eq:weq} for $w_i$ in \python{} using \fenics{} and compute the solution $u_i$ as desired. In this way, we obtain the correct solution $u_i$ at each node, but there still happens a smoothing between the nodes around the discontinuity, as we assign $u_i$ to a $H^1$ function space using Lagrange basis functions.

\subsection{Reconstruction of $\theta$}
We compute the true angle $\theta$ as the argument of $\nabla u_1$ as highlighted in equation \eqref{eq:thetau1}. Using the true angle as boundary condition for the boundary value problem \eqref{eq:thetaBVP} we reconstruct $\theta$ by solving the problem numerically. This is repeated for the two conductivities as in figure \ref{fig:TrueSig}, all different boundaries of control as in figure \ref{fig:GammaSize} and the continuous and discontinuous boundary conditions as in figure \ref{fig:contBC} and \ref{fig:discontBC}. The relative errors are shown in table \ref{tab:RelErrSize} and table \ref{tab:thetaCosSin}. 

\begin{table}[ht!]
    \centering
    \caption{Relative $L^2$ errors when using the continuous boundary conditions ($u_i^c\vert_{\partial \Omega}$) and the discontinuous boundary conditions ($u_i^d\vert_{\partial \Omega}$).}
    \begin{tabular}{c c||c c|c c|c c}
        && \multicolumn{2}{c|}{$\Gamma_{\text{large}}$} & \multicolumn{2}{c|}{$\Gamma_{\text{medium}}$} & \multicolumn{2}{c}{$\Gamma_{\text{small}}$}\\
         && $u_i^c\vert_{\partial \Omega}$ & $u_i^d\vert_{\partial \Omega}$& $u_i^c\vert_{\partial \Omega}$ & $u_i^d\vert_{\partial \Omega}$& $u_i^c\vert_{\partial \Omega}$ & $u_i^d\vert_{\partial \Omega}$\\\hline \hline
        \multirow{2}{*}{Min det$(\mathbf{H})$} & case 1 & $1 \cdot 10^{-6}$ & $3 \cdot 10^{-6}$ & $8 \cdot 10^{-10}$ & $1 \cdot 10^{-9}$ & $8 \cdot 10^{-14}$ & $8 \cdot 10^{-14}$\\
        & case 2 & $1 \cdot 10^{-6}$ & $3 \cdot 10^{-6}$ & $9 \cdot 10^{-10}$ & $2 \cdot 10^{-9}$ & $8 \cdot 10^{-14}$ & $7 \cdot 10^{-14}$\\\hline
        \multirow{2}{*}{Rel. $L^2$ error $\theta$} & case 1 & 1.62\% & 0.74\% & 1.19\% & 6.90\% & - & -\\
         & case 2 & 1.67\% & 0.75\% & 1.20\% & 7.11\% & - & -\\ \hline
        \multirow{2}{*}{Rel. $L^2$ error $\sigma$} & case 1 & 15.7\% & 15.6\% & 40.1\% & 39.9\% & 56.5\% & 56.3\%\\
        &case 2 & 15.0\% & 15.0\% & 40.0\% & 39.9\% & 56.9\% & 56.4\%\\\hline
    \end{tabular}
    \label{tab:RelErrSize}
\end{table}

\begin{table}[ht!]
    \centering
    \caption{Relative $L^2$ errors of $(\cos(2\theta),\sin(2\theta))$ for $\Gamma_{\text{small}}$.}
    \begin{tabular}{c||c c}
         & Continuous BC & Discontinuous BC\\ \hline \hline
         case 1 & (3.47\%, 3.61\%) & (5.34\%, 5.44\%)\\
         case 2 & (3.58\%, 3.56\%) & (5.40\%, 5.38\%)\\ \hline
    \end{tabular}
    \label{tab:thetaCosSin}
\end{table}

\begin{figure}[ht!]
    \centering
    \begin{minipage}[t]{0.03\textwidth}
        \begin{turn}{90}\hspace{6mm}Continuous BC\end{turn}
    \end{minipage}%
    \begin{minipage}[t]{0.33\textwidth}
        \centering
        \includegraphics[trim={0cm 0cm 4cm 0cm},clip,width=0.76\linewidth]{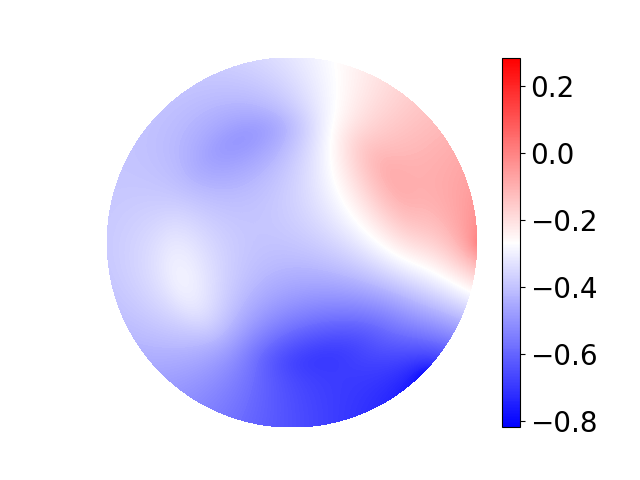}
    \end{minipage}%
    \begin{minipage}[t]{0.33\textwidth}
        \centering
        \includegraphics[width=\linewidth]{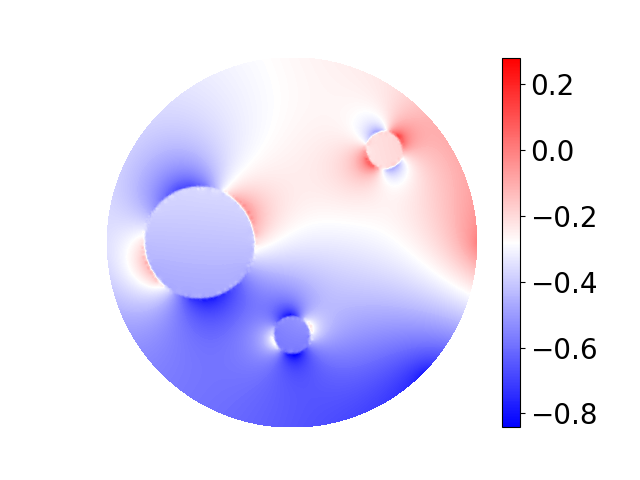}
    \end{minipage}
    \begin{minipage}[t]{0.03\textwidth}
        \begin{turn}{90}\hspace{3mm}Discontinuous BC\end{turn}
    \end{minipage}%
    \begin{minipage}[t]{0.33\textwidth}
        \centering
        \includegraphics[trim={0cm 0cm 4cm 0cm},clip,width=0.76\linewidth]{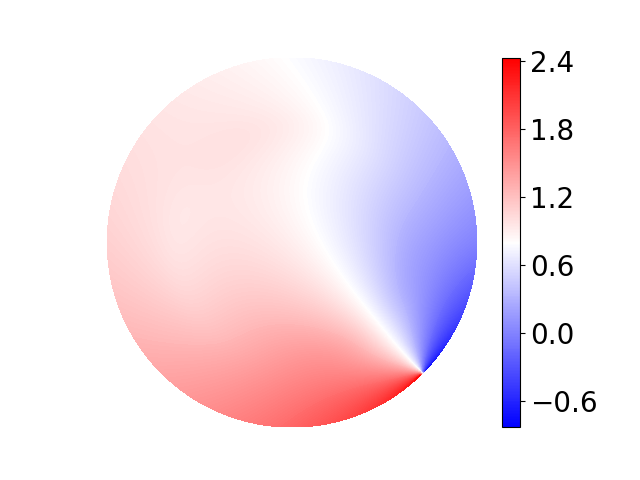}
        \caption*{\hspace{1cm}Case 1}
    \end{minipage}%
    \begin{minipage}[t]{0.33\textwidth}
        \centering
        \includegraphics[width=\linewidth]{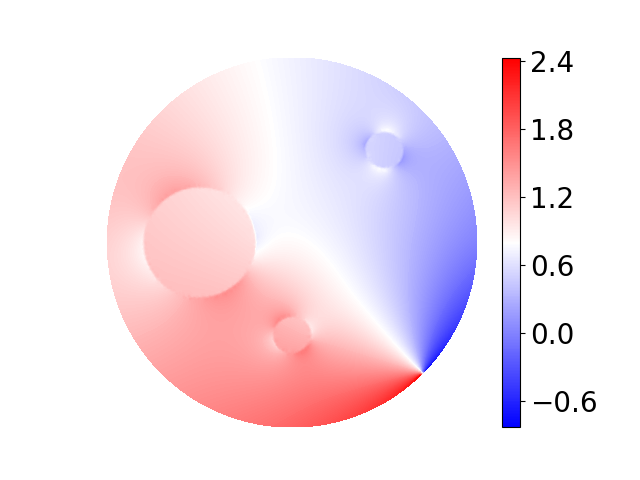}
        \caption*{Case 2}
    \end{minipage}
        \caption{True $\theta$ for the different boundary conditions and different conductivities and having control over $\Gamma_{\text{large}}$. }
    \label{fig:thetaLar}
\end{figure}

\begin{figure}[ht!]
    \centering
    \begin{minipage}[t]{0.03\textwidth}
        \begin{turn}{90}\hspace{6mm}Continuous BC\end{turn}
    \end{minipage}%
    \begin{minipage}[t]{0.33\textwidth}
        \centering
        \includegraphics[trim={0cm 0cm 4cm 0cm},clip,width=0.76\linewidth]{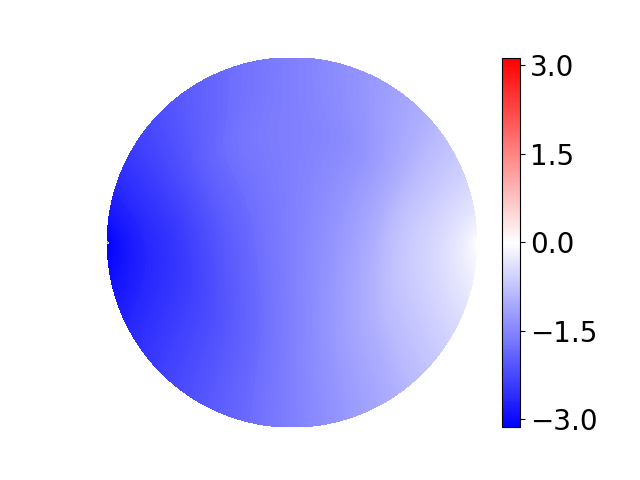}
    \end{minipage}%
    \begin{minipage}[t]{0.33\textwidth}
        \centering
        \includegraphics[width=\linewidth]{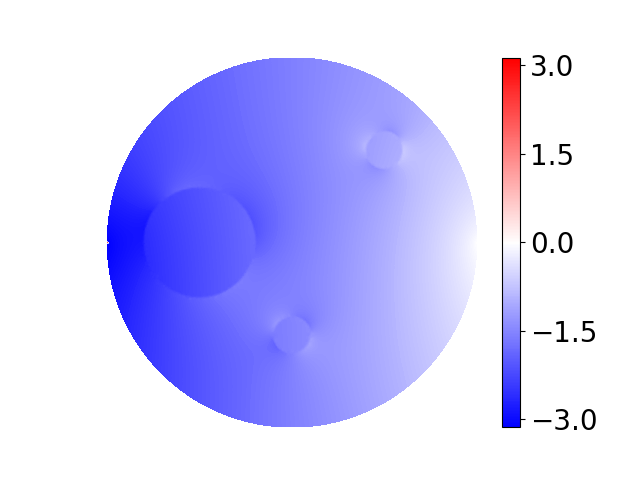}
    \end{minipage}
    \begin{minipage}[t]{0.03\textwidth}
        \begin{turn}{90}\hspace{3mm}Discontinuous BC\end{turn}
    \end{minipage}%
    \begin{minipage}[t]{0.33\textwidth}
        \centering
        \includegraphics[trim={0cm 0cm 4cm 0cm},clip,width=0.76\linewidth]{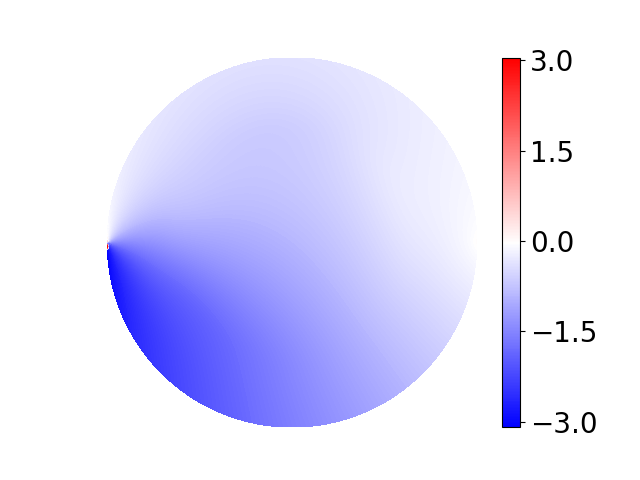}
        \caption*{\hspace{1cm}Case 1}
    \end{minipage}%
    \begin{minipage}[t]{0.33\textwidth}
        \centering
        \includegraphics[width=\linewidth]{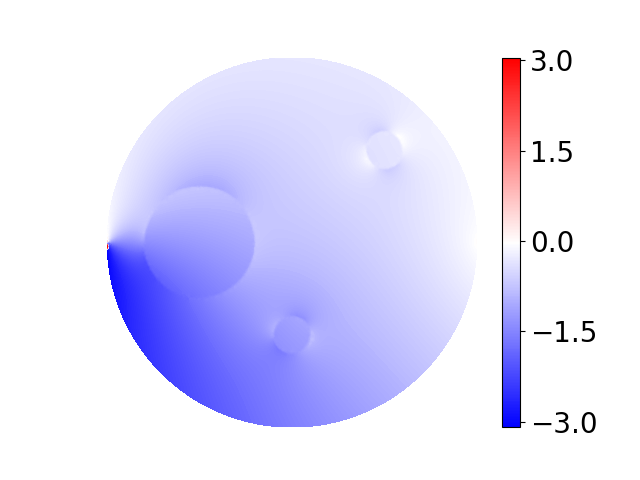}
        \caption*{Case 2}
    \end{minipage}
    \caption{True $\theta$ for the different boundary conditions and different conductivities and having control over $\Gamma_{\text{medium}}$.}
    \label{fig:thetaMed}
\end{figure}

Even though the errors range up to $7\%$ the reconstructions with control over $\Gamma_{\text{large}}$ or $\Gamma_{\text{medium}}$ can barely be distinguished visually from the true $\theta$. The only difference appears through minor artifacts along the part of the boundary that cannot be controlled. For that purpose in these cases, we only focus on the true expressions for $\theta$. These are illustrated for $\Gamma_{\text{large}}$ in figure \ref{fig:thetaLar} and for $\Gamma_{\text{medium}}$ in figure \ref{fig:thetaMed}. From the figures we see that the discontinuities of $u_i^d$ are reflected in $\theta$ as well, which follows directly from the definition of $\theta$. Furthermore, $\theta$ differs for the two cases of $\sigma$: As the circular feature in $\sigma_{\text{case 1}}$ has a smooth edge its contours can barely be seen in the expression for $\theta$. However, as there appear piecewise constant features in $\sigma_{\text{case 2}}$ the edges are clearly reflected in the expression for $\theta$ as well. \par
\begin{figure}[!ht]
    \centering
    \begin{minipage}[t]{0.03\textwidth}
        \begin{turn}{90}\hspace{13mm}Continuous BC\end{turn}
    \end{minipage}%
    \begin{minipage}[t]{0.5\textwidth}
        \centering
        \includegraphics[width=\linewidth]{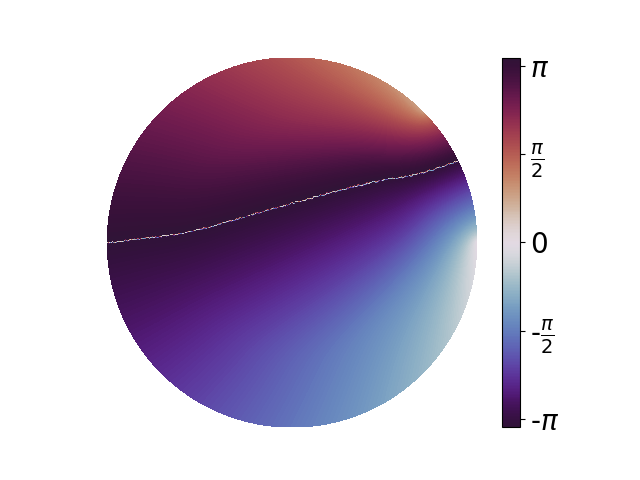}
    \end{minipage}%
    \begin{minipage}[t]{0.4\textwidth}
        \centering
        \includegraphics[width=\linewidth]{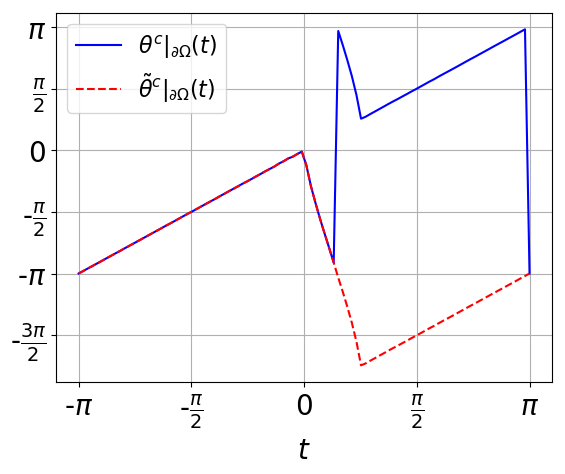}
    \end{minipage}
    \begin{minipage}[t]{0.03\textwidth}
        \begin{turn}{90}\hspace{10mm}Discontinuous BC\end{turn}
    \end{minipage}%
    \begin{minipage}[t]{0.5\textwidth}
        \centering
        \includegraphics[width=\linewidth]{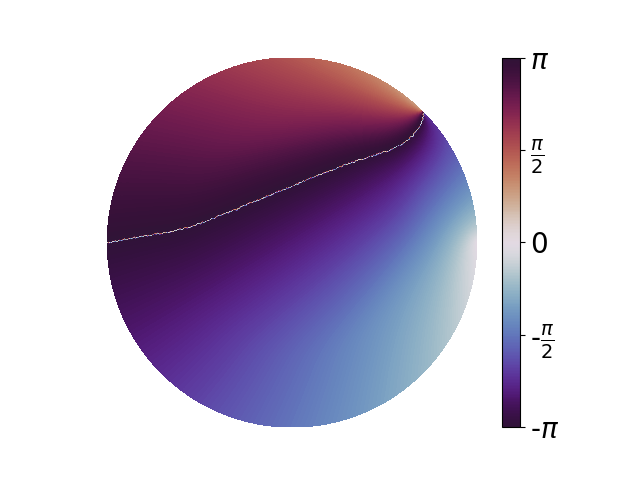}
        \caption*{True $\theta$}
    \end{minipage}%
    \begin{minipage}[t]{0.4\textwidth}
        \centering
        \includegraphics[width=\linewidth]{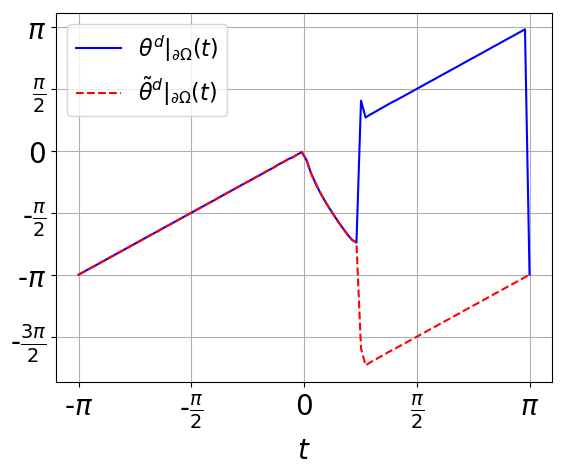}
        \caption*{Comparing $\theta$ to a modified version $\tilde{\theta}$ along $\partial \Omega$}
    \end{minipage}
    \caption{True expression for $\theta$ assigned to a smooth function space and using continuous and discontinuous boundary conditions. We consider the conductivity $\sigma_{\text{case 1}}$ and have control over $\Gamma_{\text{small}}$. The left part shows $\theta$ along the boundary together with modified versions $\tilde{\theta}$ defined in \eqref{eq:tthetac} and \eqref{eq:tthetad}. }
    \label{fig:thetaSma}
\end{figure}

 
For some choices of the boundary conditions and $\Gamma$, the true angle $\theta$ changes values from $-\pi$ to $\pi$ throughout $\Omega$. This is here the case for both continuous and discontinuous boundary conditions and when $\Gamma$ has the size of $\Gamma_{\text{small}}$. This behavior causes a curve of transitions. Along this curve the expression transitions through all values from $-\pi$ to $\pi$, instead of leaving a discontinuity, which would adhere to the periodic nature of the codomain. This is illustrated in the left part of figure \ref{fig:thetaSma} and a periodic color map is used to highlight the transition curve. A similar phenomenon was observed and addressed in \cite{bjoernArxiv} and we use the same approach to address this issue. By the discussion in section \ref{sec:ChoiceTtheta} the direction of $\nabla u_1$ corresponds to the direction of the unit normal $\nu$, its opposite $-\nu$, or the zero vector along the boundary $\partial \Omega \backslash \Gamma_{\text{small}}$. If we investigate $\theta$ along the boundary in the right part of figure \ref{fig:thetaSma}, we see that along $\partial \Omega \backslash \Gamma_{\text{small}}$ $\theta$ is a linear increasing function, as it corresponds to the angle between $\nu$ and the $x_1$-axis. The only deviations happen for $t \in [0,\frac{\pi}{4}]$, which is the boundary of control $\Gamma_{\text{small}}$. By this behavior of $\theta$ along $\partial \Omega \backslash \Gamma_{\text{small}}$ there happens a jump from $\pi$ to $-\pi$ at $t=\pi$. An additional jump is induced by the boundary condition: For the continuous boundary condition there happens a jump at $t=\frac{\pi}{8}$ and for the discontinuous boundary condition there happens a jump at $t=\frac{\pi}{4}$. The smoothed discontinuities are a problem when using the true angle $\theta$ as a boundary condition in \eqref{eq:thetaBVP}, therefore we define a modified version $\tilde{\theta}$ to be used as a boundary condition. As $\theta$ only appears as an input to the cosine and sine-functions in the reconstruction formula in \eqref{eq:sigmaBVP}, we can add and subtract multiples of $2\pi$ without changing the reconstruction. Therefore, we subtract $2\pi$ in the interval between the discontinuities to extend $\theta$ to a more continuous function along the boundary. For continuous boundary conditions $\tilde{\theta}^c$ is defined as 
\begin{equation}\label{eq:tthetac}
    \tilde{\theta}^c(t) = \begin{cases}
        \theta^c(t) - 2\pi & t \in \left[\frac{\pi}{8},\pi \right],\\
        \theta^c(t) & \text{otherwise}.
    \end{cases}
\end{equation}
And for the discontinuous boundary condition:
\begin{equation}\label{eq:tthetad}
    \tilde{\theta}^d(t) = \begin{cases}
        \theta^d(t) - 2\pi & t \in \left[\frac{\pi}{4},\pi \right],\\
        \theta^d(t) & \text{otherwise}.
    \end{cases}
\end{equation}
These functions are illustrated in the right part of figure \ref{fig:thetaSma} and used as a boundary condition when solving the boundary value problem \eqref{eq:thetaBVP}. 

\begin{figure}[ht!]
    \centering
    \begin{minipage}[t]{0.03\textwidth}
        \begin{turn}{90}\hspace{6mm}Continuous BC\end{turn}
    \end{minipage}%
    \begin{minipage}[t]{0.33\textwidth}
        \centering
        \includegraphics[trim={0cm 0cm 4cm 0cm},clip,width=0.76\linewidth]{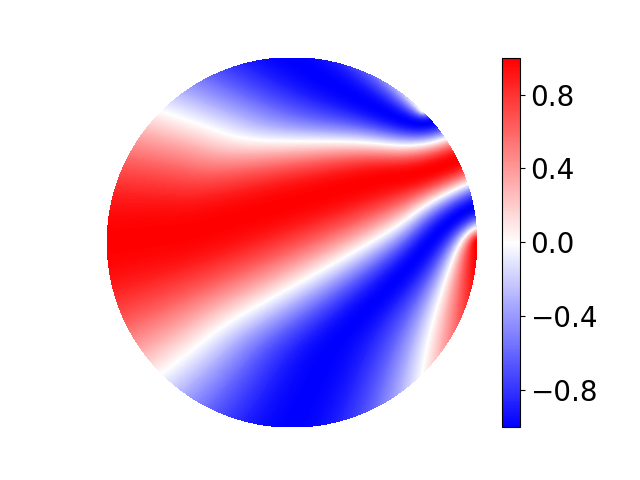}
    \end{minipage}%
    \begin{minipage}[t]{0.33\textwidth}
        \centering
        \includegraphics[width=\linewidth]{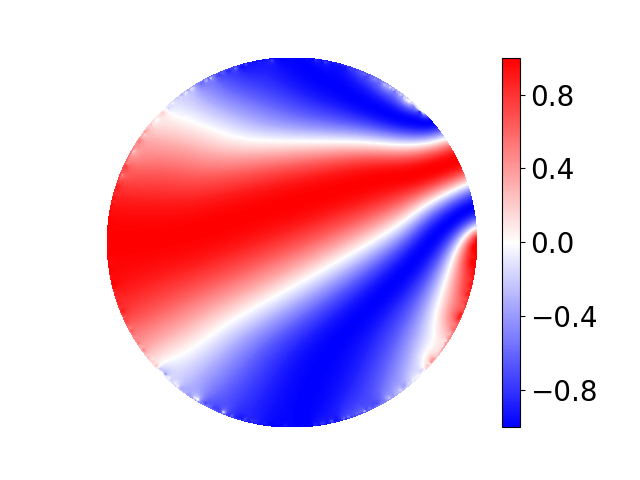}
    \end{minipage}
    \begin{minipage}[t]{0.03\textwidth}
        \begin{turn}{90}\hspace{3mm}Discontinuous BC\end{turn}
    \end{minipage}%
    \begin{minipage}[t]{0.33\textwidth}
        \centering
        \includegraphics[trim={0cm 0cm 4cm 0cm},clip,width=0.76\linewidth]{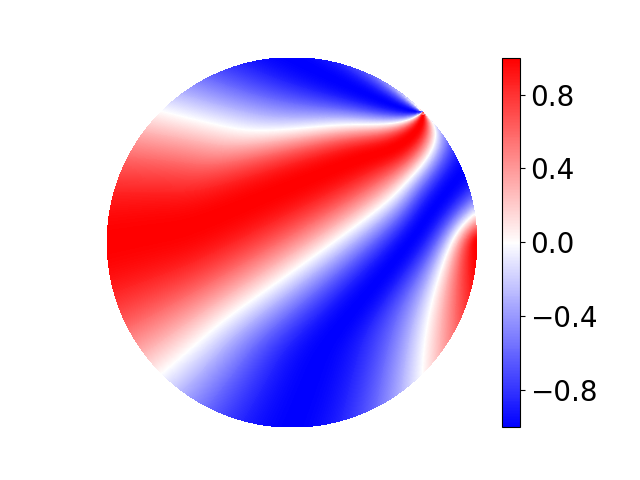}
        \caption*{\hspace{1cm}True $\cos(2\theta)$}
    \end{minipage}%
    \begin{minipage}[t]{0.33\textwidth}
        \centering
        \includegraphics[width=\linewidth]{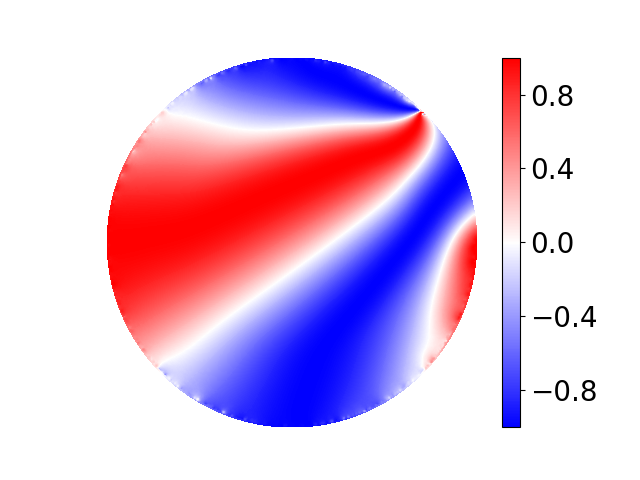}
        \caption*{Reconstructed $\cos(2\theta)$}
    \end{minipage}
    \caption{Reconstructions of $\cos(2\theta)$ compared to the true expression for the different boundary conditions. We consider the conductivity $\sigma_{\text{case 1}}$ and have control over $\Gamma_{\text{small}}$.}
    \label{fig:thetaCos}
\end{figure}

For these reconstructions $\theta$ it does no longer make sense to compare them to the true angles, so instead we compare the reconstructed $\cos(2\theta)$ to the true expression in figure \ref{fig:thetaCos}, as this is the way $\theta$ appears in the reconstruction formula for $\sigma$ \eqref{eq:sigmaBVP}. The reconstruction errors are shown in table \ref{tab:thetaCosSin}. From the figure, we see that there appear artifacts along the whole boundary $\partial \Omega \backslash \Gamma_{\text{small}}$. This is because this part of the boundary is difficult to control from the small boundary $\Gamma_{\text{small}}$ so that the Jacobian constraint almost is violated close to this part of the boundary. This is in accordance with remark \ref{rem:vio}, as the Jacobian constraint is violated on this part of the boundary. This is seen from figure \ref{fig:Jac}, as the values of $\text{det}(\mathbf{H})$ are very small close to the boundary $\partial \Omega \backslash \Gamma_{\text{small}}$. And from table \ref{tab:RelErrSize} we see the small values are of order $10^{-14}$. These artifacts were not as visible for $\Gamma_{\text{large}}$ and $\Gamma_{\text{medium}}$, as the smallest values of $\text{det}(\mathbf{H})$ were larger than $3 \cdot 10^{-6}$ and $9 \cdot 10^{-10}$ respectively. 

\subsection{Reconstruction of $\sigma$}
\begin{figure}[ht!]
    \centering
    \begin{minipage}[t]{0.33\textwidth}
        \centering
        \includegraphics[width=\linewidth]{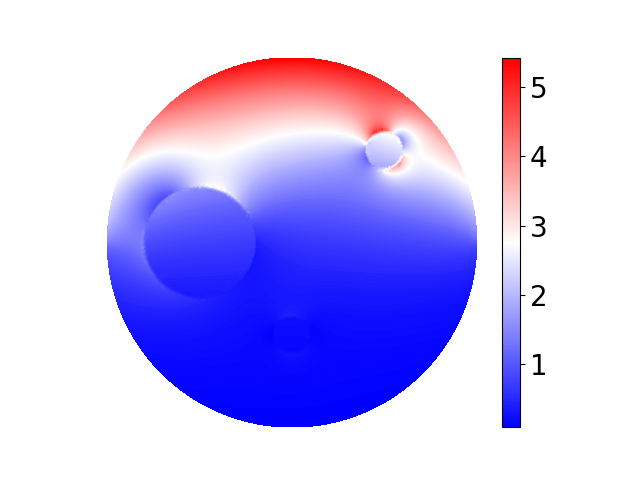}
        \caption*{$H_{11}$}
    \end{minipage}%
    \begin{minipage}[t]{0.33\textwidth}
        \centering
        \includegraphics[width=\linewidth]{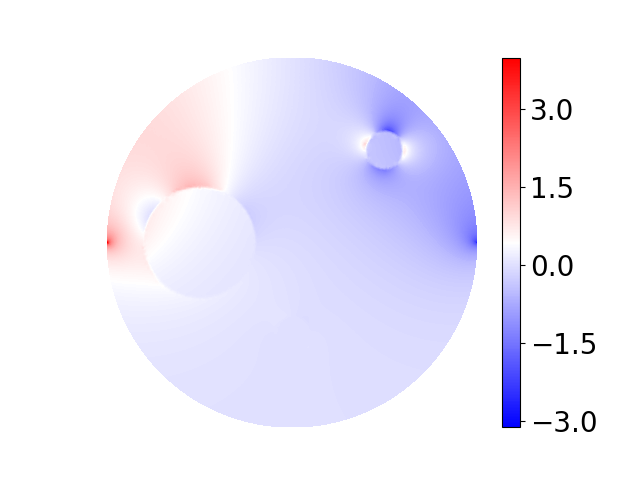}
        \caption*{$H_{12}$}
    \end{minipage}%
    \begin{minipage}[t]{0.33\textwidth}
        \centering
        \includegraphics[width=\linewidth]{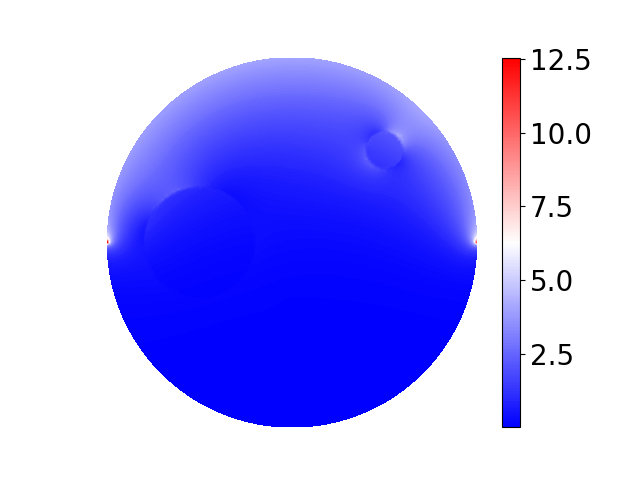}
        \caption*{$H_{22}$}
    \end{minipage}
     \caption{The power densities used for reconstruction of $\sigma_{\text{case 2}}$ and having control over $\Gamma_{\mathrm{medium}}$.}
    \label{fig:Hijs}
\end{figure}

\begin{figure}[ht!]
    \centering
    \begin{minipage}[t]{0.03\textwidth}
        \begin{turn}{90}\hspace{20mm}True $\sigma$\end{turn}
    \end{minipage}%
    \begin{minipage}[t]{0.45\textwidth}
        \centering
        \includegraphics[width=\linewidth]{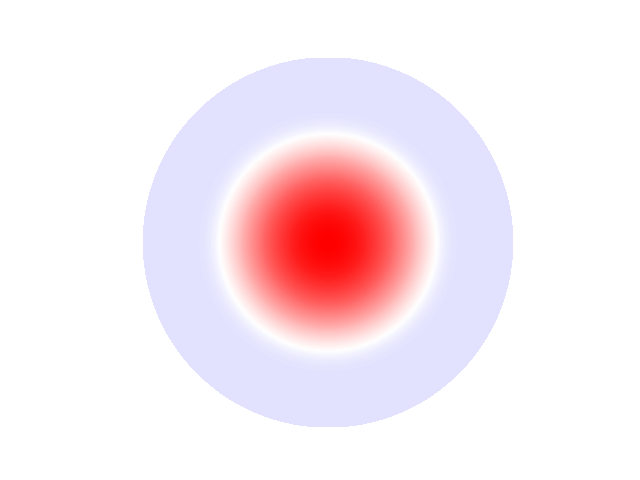}
    \end{minipage}%
    \begin{minipage}[t]{0.45\textwidth}
        \centering
        \includegraphics[width=\linewidth]{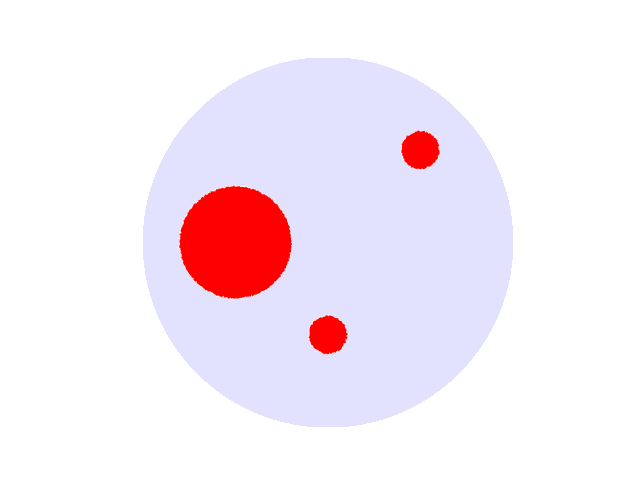}
    \end{minipage}
    \begin{minipage}[t]{0.03\textwidth}
        \begin{turn}{90}\hspace{20mm}$\Gamma_{\text{large}}$\end{turn}
    \end{minipage}%
    \begin{minipage}[t]{0.45\textwidth}
        \centering
        \includegraphics[width=\linewidth]{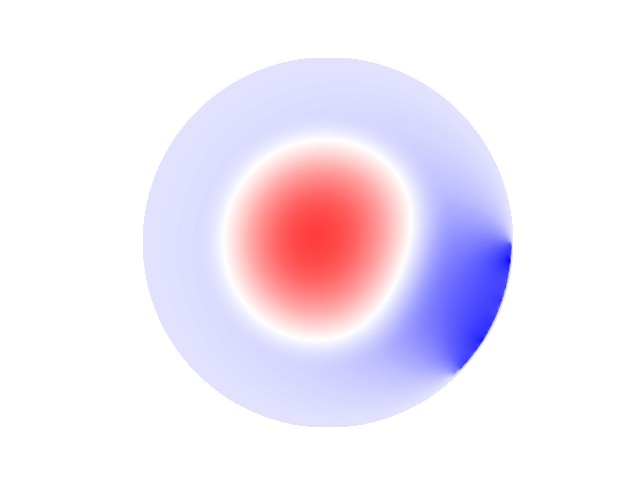}
    \end{minipage}%
    \begin{minipage}[t]{0.45\textwidth}
        \centering
        \includegraphics[width=\linewidth]{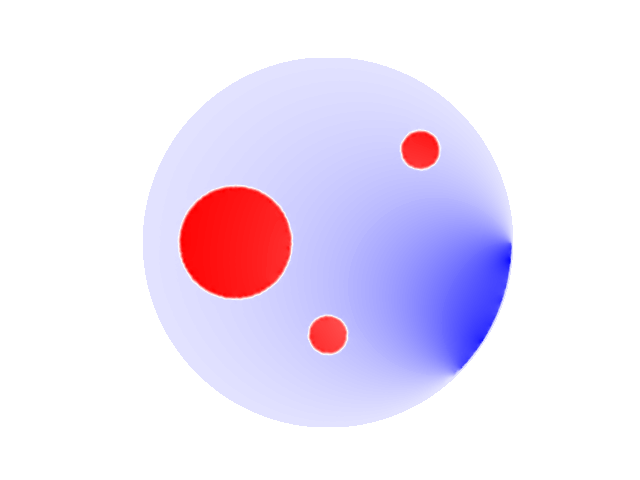}
    \end{minipage}
    \begin{minipage}[t]{0.03\textwidth}
        \begin{turn}{90}\hspace{20mm}$\Gamma_{\text{medium}}$\end{turn}
    \end{minipage}%
    \begin{minipage}[t]{0.45\textwidth}
        \centering
        \includegraphics[width=\linewidth]{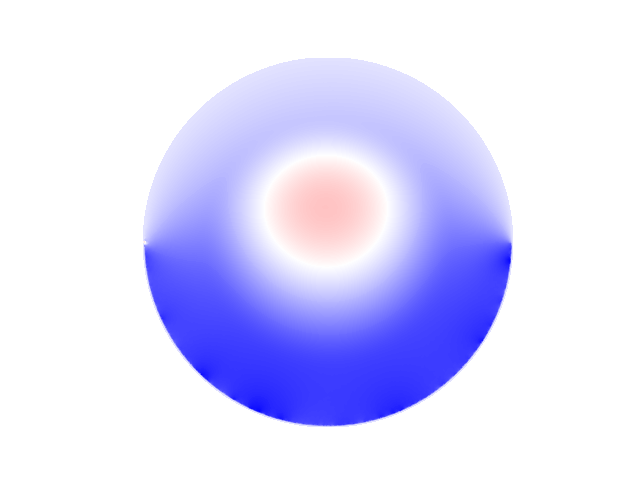}
    \end{minipage}%
    \begin{minipage}[t]{0.45\textwidth}
        \centering
        \includegraphics[width=\linewidth]{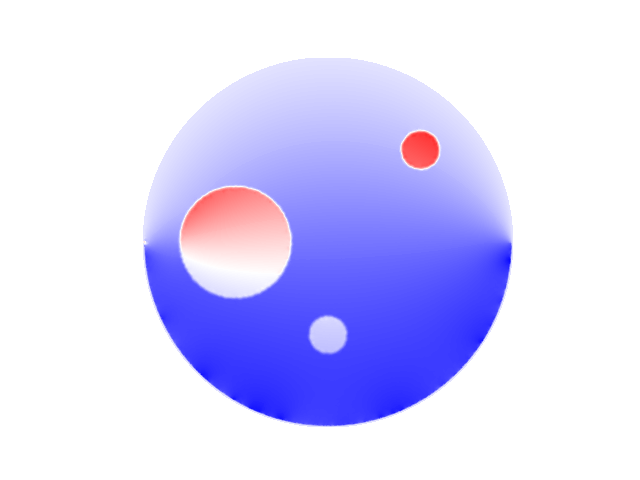}
    \end{minipage}
    \begin{minipage}[t]{0.03\textwidth}
        \begin{turn}{90}\hspace{24mm}$\Gamma_{\text{small}}$\end{turn}
    \end{minipage}%
    \begin{minipage}[t]{0.03\textwidth}
    {\color{white} bla}
    \end{minipage}%
    \begin{minipage}[t]{0.44\textwidth}
        \centering
        \includegraphics[width=0.9\linewidth]{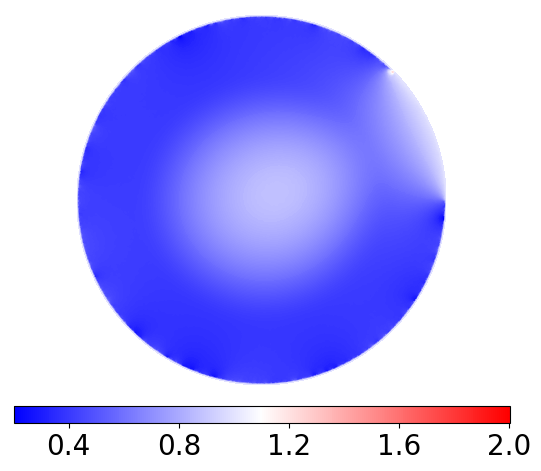}
        \caption*{Case 1}
    \end{minipage}%
    \begin{minipage}[t]{0.44\textwidth}
        \centering
        \includegraphics[width=0.9\linewidth]{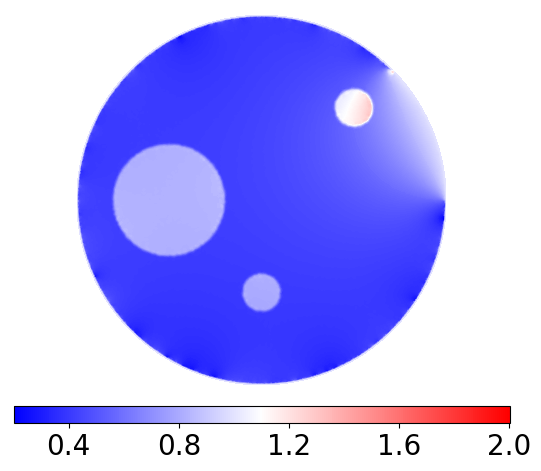}
        \caption*{Case 2}
    \end{minipage}
    \caption{Reconstructions of $\sigma$ as in test case 1 in the left column and as in test case 2 in the right column for varying sizes of $\Gamma$. The discontinuous boundary conditions are used. }
    \label{fig:sigmaRecs}
\end{figure}

Using the reconstructions of $\theta$ we compute $\sigma$ by solving the boundary value problem \eqref{eq:sigmaBVP}. The reconstructions of $\sigma$ using the continuous boundary conditions are shown in figure \ref{fig:sigmaRecs}. As seen from the relative errors in table \ref{tab:RelErrSize} there is no significant difference in the quality of the reconstruction when using the continuous or discontinuous boundary conditions. It is therefore impossible to distinguish the reconstructions visually, so we only show the reconstructions using the continuous boundary conditions. From the figure we see that the quality of the reconstructions is highly affected by the size of the boundary of control $\Gamma$: The larger the boundary of control the better the reconstruction. We note that since the minimum values of the determinant of $\mathbf{H}$ are close to machine precision for $\Gamma_{\text{small}}$, as can be seen from the first to rows in \ref{tab:RelErrSize}, the numerical results for this size of $\Gamma$ are not very reliable. We still include this example to show the limitations of the method. We see from figure \ref{fig:sigmaRecs} that for $\Gamma_{\text{large}}$ the features in the reconstructions are still visible in almost the same intensity as for the true $\sigma$, only the shape of the circular feature in $\sigma_{\text{case 1}}$ is changed in the direction of $\partial \Omega \backslash \Gamma_{\text{large}}$. With decreasing size of $\Gamma$ the intensity of the features is decreasing as well, so that for $\Gamma_{\text{small}}$ the features have intensity close to 1 like the background of the true $\sigma$. Also in a large neighborhood of the boundary $\partial \Omega \backslash \Gamma$, the reconstruction has intensity close to 0, so that for $\Gamma_{\text{small}}$ the reconstruction is dominated by intensity 0. When comparing performance for the two different conductivities $\sigma$ we see that the reconstructions for the piecewise constant conductivity $\sigma_{\text{case 2}}$ look almost better than for the smooth $\sigma_{\text{case 1}}$, as the piecewise constant edges of the three features in $\sigma_{\text{case 2}}$ are clearly visible in the reconstructions. This is due to the fact, that these edges are clearly visible in the data (see figure \ref{fig:Hijs}) and in $\theta$ (see the right parts of figure \ref{fig:thetaLar} and figure \ref{fig:thetaMed} respectively). For $\sigma_{\text{case 1}}$ the shape of the feature is deformed a little bit towards $\partial \Omega \backslash \Gamma$, as the feature has a smooth edge. However, this difference in quality is not evident from the relative errors in table \ref{tab:RelErrSize}. Another take away from the reconstructions is that as $\sigma_{\text{case 2}}$ is composed of features that are closer and further away from the boundary of control as the feature in $\sigma_{\text{case 1}}$, we can see that there is a difference in the intensity of the three features. This is especially evident for $\Gamma_{\text{medium}}$, so that the feature closest to $\Gamma_{\text{medium}}$ has intensity 1.8, which is almost the same intensity as the true $\sigma_{\text{case 2}}$. On the other hand, the feature furthest away from $\Gamma_{\text{medium}}$ has intensity 1, which is the same as the background intensity of the true $\sigma_{\text{case 2}}$.

\subsection{Reconstruction of $\sigma$ from noisy data}
We perturb the entries of the power density matrix $\mathbf{H}$ at each node with random noise:
\begin{equation*}
    \widetilde{H}_{ij} = H_{ij} + \frac{\alpha}{100} \frac{e_{ij}}{\norm{e_{ij}}_{L^2}} H_{ij},
\end{equation*}
where $\alpha$ is the noise level and $e_{ij}$ are entries in the matrix $\mathbf{E}$ that are normally distributed $e_{ij} \sim \mathcal{N}(0,1)$. We use \texttt{numpy.random.randn} to generate the elements $e_{ij}$ and fix the seed \texttt{numpy.random.seed(50)}. After generating $\widetilde{\mathbf{H}}$, we make sure that it is symmetric by computing $\frac{1}{2}(\widetilde{\mathbf{H}}+\widetilde{\mathbf{H}}^T)$. Furthermore, for the reconstruction procedure it is essential that $\widetilde{\mathbf{H}}$ is positive definite so that we use a small positive lower bound $L$ for the eigenvalues of $\widetilde{\mathbf{H}}$. This approach can be seen as a regularization method. We note that by remark \ref{rem:vio} the Jacobian constraint is violated on $\partial \Omega \backslash \Gamma$, therefore we would assume very small values of $\det(\widetilde{\mathbf{H}})$ close to this part of the boundary. However, in the approach of using a lower bound for the eigenvalues we might discard some of these values. Therefore, we choose the lower bound as small as possible in order to get a reasonable reconstruction that is not dominated by noise on these small values. \par
\begin{figure}[ht!]
    \centering
    \begin{minipage}[t]{0.45\textwidth}
        \centering
        \includegraphics[trim={0cm 1.5cm 0cm 0cm},width=\linewidth]{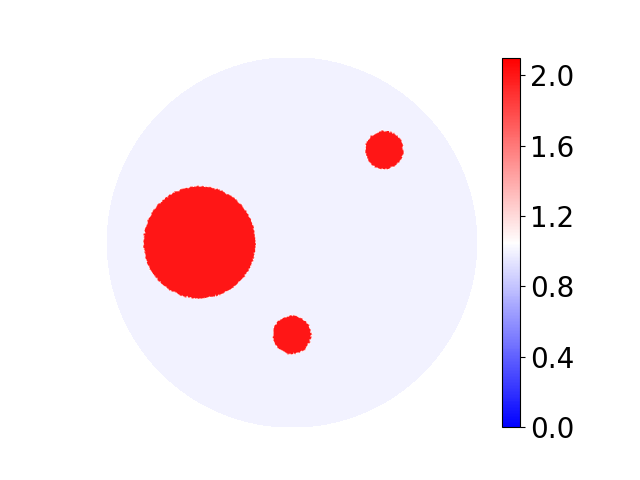}
        \caption*{True $\sigma$}
    \end{minipage}
    \begin{minipage}[t]{0.03\textwidth}
        \begin{turn}{90}\hspace{17mm}1\% Noise\end{turn}
    \end{minipage}%
    \begin{minipage}[t]{0.45\textwidth}
        \centering
        \includegraphics[trim={0cm 0cm 4cm 0cm},clip,width=0.76\linewidth]{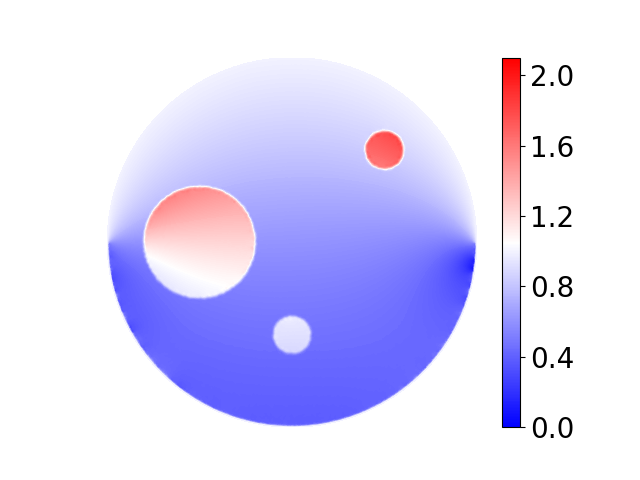}
    \end{minipage}%
    \begin{minipage}[t]{0.45\textwidth}
        \centering
        \includegraphics[trim={0cm 0cm 4cm 0cm},clip,width=0.75\linewidth]{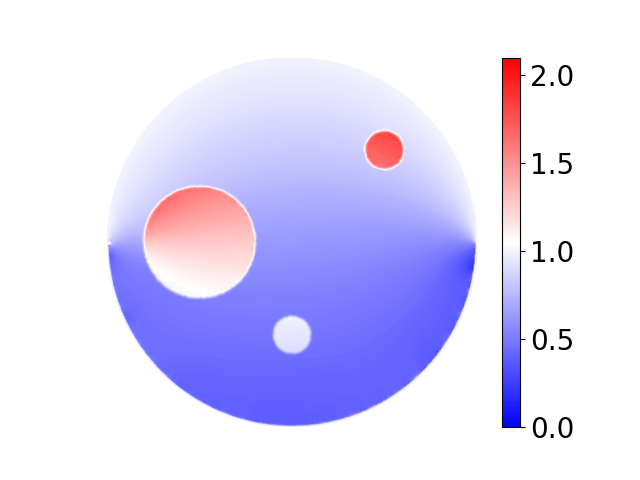}
    \end{minipage}
    \begin{minipage}[t]{0.03\textwidth}
        \begin{turn}{90}\hspace{12mm}5\% Noise\end{turn}
    \end{minipage}%
    \begin{minipage}[t]{0.45\textwidth}
        \centering
        \includegraphics[trim={0cm 0cm 4cm 0cm},clip,width=0.76\linewidth]{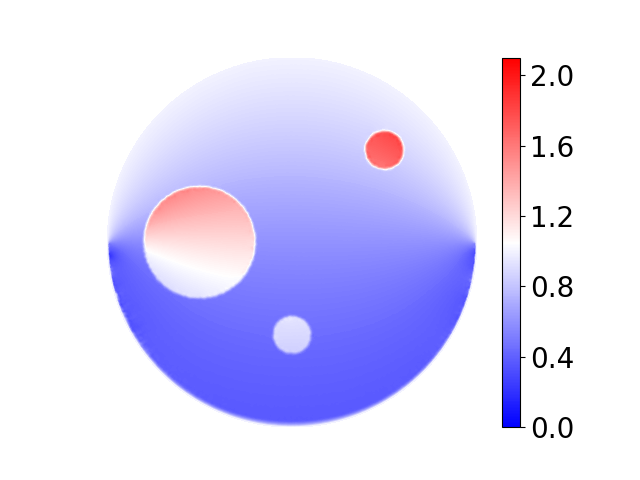}
    \end{minipage}%
    \begin{minipage}[t]{0.45\textwidth}
        \centering
        \includegraphics[trim={0cm 0cm 4cm 0cm},clip,width=0.76\linewidth]{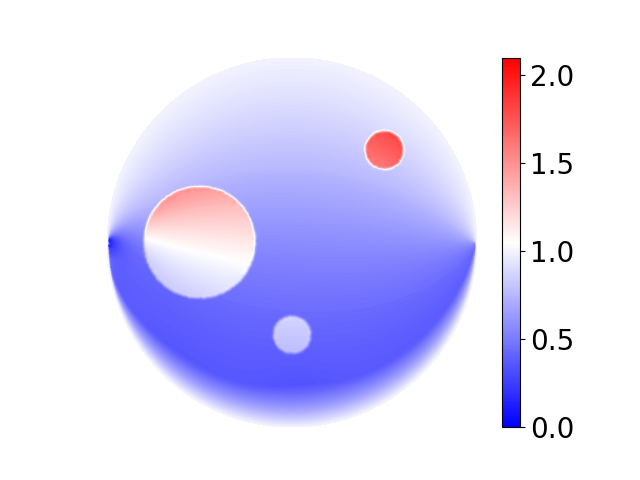}
    \end{minipage}
    \begin{minipage}[t]{0.03\textwidth}
        \begin{turn}{90}\hspace{17mm}10\% Noise\end{turn}
    \end{minipage}%
    \begin{minipage}[t]{0.45\textwidth}
        \centering
        \includegraphics[trim={0cm 0cm 4cm 0cm},clip,width=0.76\linewidth]{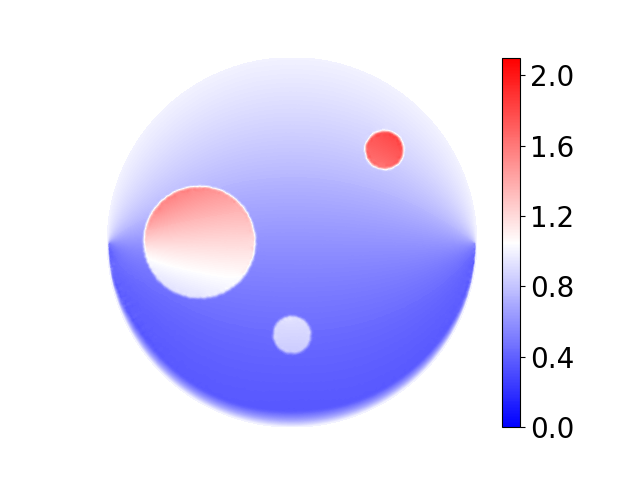}
        \caption*{\hspace{1cm}Continuous BC}
    \end{minipage}%
    \begin{minipage}[t]{0.45\textwidth}
        \centering
        \includegraphics[trim={0cm 0cm 4cm 0cm},clip,width=0.76\linewidth]{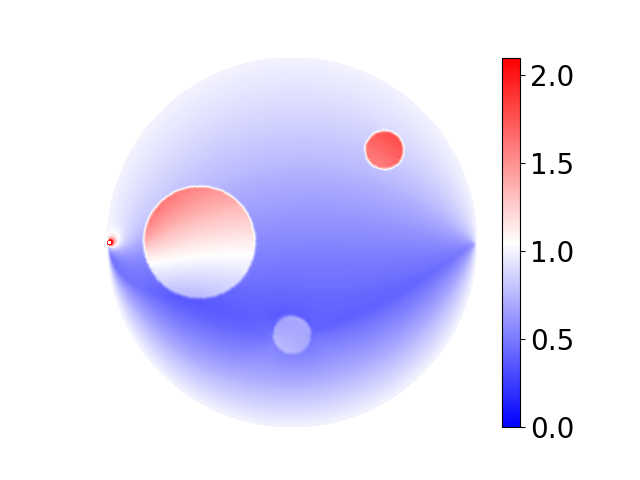}
        \caption*{\hspace{1cm}Discontinuous BC}
    \end{minipage}
    \caption{Reconstructions of $\sigma_{\text{case 2}}$ when $\mathbf{H}$ is perturbed with varying noise levels and with boundary of control $\Gamma_{\mathrm{medium}}$. To obtain a positive definite noisy matrix $\widetilde{\mathbf H}$, different lower bounds $L$ for the eigenvalues of $\widetilde{\mathbf H}$ are used. }
    \label{fig:sigmaRecsnoise}
\end{figure}

After these modifications on symmetry and positive definiteness, we use $\widetilde{\mathbf{H}}$ for reconstructing $\sigma_{\text{case 2}}$ for three different noise levels: $\alpha=1$, $\alpha=5$ and $\alpha=10$. The results are shown in Figure~\ref{fig:sigmaRecsnoise}, where we compare performance when using continuous and discontinuous boundary conditions. The lower bounds $L$ used for reconstruction are documented in Table~\ref{tab:RelErrNoise} and the relative errors of $\sigma$ are shown in the same table. The reconstructions are all of similar quality when looking at the relative errors and how well the values at the red features matches with the true $\sigma$. However, for the discontinuous boundary condition the increasing noise level results in an artifact around the discontinuity. To account for this one needs to use significantly larger lower bounds for increasing noise level. In contrast, there is a gradual rise in the lower bound for increasing noise level when using continuous boundary conditions. The high lower bound in the case of discontinuous boundary conditions induces a light belt close to $\partial \Omega \backslash \Gamma$. This belt appears, as information in this region is discarded by the lower bound.

\begin{table}[ht!]
    \centering
    \caption{Relative $L^2$ errors on $\sigma_{\text{case 2}}$ in presence of noise. To obtain a positive definite noisy matrix $\widetilde{\mathbf H}$, different lower bounds $L$ for the eigenvalues of $\widetilde{\mathbf H}$ are used. The boundary of control is $\Gamma_{\mathrm{medium}}$.}
    \begin{tabular}{c||c | c|c | c}
    &\multicolumn{2}{c|}{Continous BC} & \multicolumn{2}{c}{Discontinuous BC}\\
         & $L$ & Relative error $\sigma$ & $L$ & Relative error $\sigma$\\\hline \hline
       1\% Noise & $10^{-6}$ & 40.7\% & $10^{-6}$ & $38.9\%$\\\hline
        5\% Noise & $10^{-5}$ & 41.4\% & $10^{-3}$ & $41.7\%$\\\hline
        10\% Noise & $10^{-4}$ & 40.6\% & $2\cdot 10^{-2}$ & $38.4\%$\\\hline
    \end{tabular}
    \label{tab:RelErrNoise}
\end{table}

\clearpage

\section{Conclusions}
In this work, we have derived sufficient conditions on two boundary functions so that the corresponding solutions to the conductivity equation satisfy a non-vanishing Jacobian constraint in limited view. This approach allows for boundary functions that have discontinuities. This is relevant for Acousto-Electric Tomography and Current Density Imaging both in limited view and in full view settings, as the conditions and thus the use of discontinuous boundary functions apply for both settings. \par
We illustrated how these conditions could be used for numerical examples of reconstructing the conductivity from power density data in limited view following the approach of \cite{monard2012a}. It was evident from the numerical examples how the non-vanishing Jacobian constraint was almost violated close to the boundary that could not be controlled. This follows from the zero Dirichlet condition on this part of the boundary: Here the two corresponding solutions have both the direction of the unit normal so that the non-vanishing Jacobian constraint cannot be satisfied on this part of the boundary. Nevertheless, without noise it was possible to obtain decent reconstructions of the conductivity when the support of the boundary of control was at least half the size of the full boundary. For smaller boundaries of control the numerical results get unreliable due to small values of the Jacobian. Therefore in general it is almost impossible to add even small levels of noise while maintaining positive definiteness of the measurement matrix. To account for this, we used a small positive lower bound for the eigenvalues of the measurement matrix. This worked well in the numerical experiments even for high noise levels. However, in this approach especially values close to the part of the boundary that cannot be controlled are affected by the lower bound. When the lower bound is high this might be in contradiction with the assumption that values should be small in this region. This is especially a problem when using discontinuous boundary conditions as the lower bound needs to be chosen large when the noise level is high. \par
We mention that the proposed conditions in order to obtain solutions satisfying the non-vanishing Jacobian constraint are not optimal. Especially for the case of discontinuous boundary functions there is a possibility that \ref{thm:Main} (b) can be relaxed, as we are aware of functions that are more general than allowed here as indicated in Remark \ref{rem:Discont}.

\paragraph*{Acknowledgements}
The authors would like to express their sincere gratitude to our Nuutti Hyv\"onen for suggesting how to solve the boundary value problem \eqref{bvpui} numerically with discontinuous boundary conditions, and for providing the \matlab{} code to solve the boundary value problem \eqref{eq:LapBVP} semi-analytically as discussed in section \ref{sec:DBC}. M.S.\ was partly supported by the Academy of Finland (Centre of Excellence in Inverse Modelling and Imaging, grant 284715) and by the European Research Council under Horizon 2020 (ERC CoG 770924).

\clearpage
\newpage
\printbibliography

Department of Mathematics and Statistics, University of Jyväskylä, Finland\\
\indent \textit{E-mail address: } \texttt{mikko.j.salo@jyu.fi}\\
\indent \textit{E-mail address: } \texttt{hjordis.a.schluter@jyu.fi}

\end{document}